\documentclass[12pt,a4paper]{amsart}
\usepackage[french]{babel}
\usepackage[latin1]{inputenc}
\usepackage{fancyhdr}

\usepackage{amsmath,amssymb,amsfonts}

\newcommand{\sm}[1]{{\scriptstyle #1}}
\newcommand{\N}{\mathbb{N}}
\newcommand{\Z}{\mathbb{Z}}

\newcommand{\R}{\mathbb{R}}
\newcommand{\C}{\mathbb{C}}
\newcommand{\Hb}{\mathbb{H}}
\newcommand{\Pb}{\mathbb{P}}
\newcommand{\too}{\longrightarrow}

\newcommand{\eq}{\Leftrightarrow}

\newcommand{\eps}{\varepsilon}
\newcommand{\fy}{\varphi}
\newcommand{\ges}{\geqslant}
\newcommand{\les}{\leqslant}
\newcommand{\moins}{\smallsetminus}
\newcommand{\de}{\partial}
\newcommand{\dbar}{\overline{\partial}}
\newcommand{\ddbar}{\partial\overline{\partial}}
\newcommand{\ec}{\mathcal{E}}
\newcommand{\lc}{\mathcal{L}}
\newcommand{\sca}{\mathcal{S}}
\newcommand{\hc}{\mathcal{H}}

\newcommand{\oc}{\mathcal{O}}
\newcommand{\ocbar}{\overline{\mathcal{O}}}
\newcommand{\ombar}[1]{\overline{\Omega}^{#1}}
\newcommand{\fc}{\mathcal{F}}
\newcommand{\gc}{\mathcal{G}}
\newcommand{\bc}{\mathcal{B}}
\newcommand{\dc}{\mathcal{D}}
\newcommand{\xc}{\mathcal{X}}

\newcommand{\Ub}{\mathbf{U}}

\newcommand{\bra}{\left<}
\newcommand{\ket}{\right>}
\newcommand{\dcech}{\check{\delta}}
\newcommand{\dcechg}{\delta\!\!\!\check{\delta}}
\newcommand{\alphab}{\overline{\alpha}}
\newcommand{\betab}{\overline{\beta}}
\newcommand{\gammab}{\overline{\gamma}}

\DeclareMathOperator{\im}{im}

\DeclareMathOperator{\id}{id}
\DeclareMathOperator{\Id}{Id}

\DeclareMathOperator{\pr}{pr}

\newtheorem{thm}{Th\'eor\`eme}
\newtheorem{lemma}[thm]{Lemme}
\newtheorem{cor}[thm]{Corollaire}
\newtheorem{prop}[thm]{Proposition}

\numberwithin{thm}{section}

\theoremstyle{definition}
\newtheorem*{defn}{D\'efinition}

\theoremstyle{remark}
\newtheorem*{rem}{Remarque}
\newtheorem*{rems}{Remarques}
\newtheorem*{notation}{Notation}

\begin{document}

\baselineskip15pt
\pagestyle{fancy}
\renewcommand{\headrulewidth}{0.4pt} 
\renewcommand{\footrulewidth}{0pt}
\lfoot{} \cfoot{} \rfoot{} 
\lhead[\thepage / ]{\textit{Autour de la cohomologie de Bott-Chern}} \chead{} 
\rhead[\textit{Autour de la cohomologie de Bott-Chern}]{ / \thepage}

\title{Autour de la cohomologie de Bott-Chern}
\author{Michel Schweitzer}
\maketitle

\tableofcontents

\pagebreak

\pagestyle{fancy}
\renewcommand{\headrulewidth}{0.4pt} 
\renewcommand{\footrulewidth}{0pt}
\lfoot{} \cfoot{} \rfoot{} 
\lhead[\thepage / Cohomologie de BC comme espaces de formes]{\textit{Autour de la cohomologie de Bott-Chern}} \chead{} 
\rhead[\textit{Autour de la cohomologie de Bott-Chern}]{Cohomologie de BC comme espaces de formes / \thepage}

\part*{Première partie : Groupes de cohomologie de Bott-Chern comme espaces de formes}

\section{Définition des groupes de Bott-Chern}

On commence par rappeler la définition des groupes de cohomologie de Bott-Chern usuelle sur une variété complexe lisse, définis comme espaces de formes ([Dem93]).

\subsection{Invariants cohomologiques d'une variété complexe}$\;$ 

Soit $X$ une variété analytique complexe. On considère, pour $k\in\{0,\ldots,2n=\dim_{\R}X\}$, l'espace $\ec^k(X)$ des formes différentielles de degré $k$ à valeurs complexes sur $X$. Celui-ci admet la décomposition
$$\ec^k(X)=\bigoplus_{\substack{0\les p,q\les n \\ p+q=k}}\ec^{p,q}(X)$$
où $\ec^{p,q}(X)$ désigne l'espace des formes de type $(p,q)$ sur $X$.\\ 
La différentielle $d:\ec^k(X)\to \ec^{k+1}(X)$ se décompose sous la forme $d=\de+\dbar$ avec
$$\de:\ec^{p,q}(X)\to \ec^{p+1,q}(X),\qquad \dbar:\ec^{p,q}(X)\to \ec^{p,q+1}(X).$$

Les invariants cohomologiques traditionnels sur $X$ sont les groupes de cohomologie de De Rham et de Dolbeault, définis par
$$H^k_{DR}(X,\C)=\frac{\ker(d:\ec^k(X)\to\ec^{k+1}(X))}{\im(d:\ec^{k-1}(X)\to\ec^k(X))}$$
$$H^{p,q}_{\dbar}(X,\C)=\frac{\ker(\dbar:\ec^{p,q}(X)\to\ec^{p,q+1}(X))}{\im(\dbar:\ec^{p,q-1}(X)\to\ec^{p,q}(X))}.$$

\begin{defn} Le groupe de cohomologie de Bott-Chern de bidegré $(p,q)$ sur une variété complexe $X$ est l'espace des $(p,q)$-formes $d$-fermées quotienté par l'espace des formes $\ddbar$-exactes :
$$H^{p,q}_{BC}(X,\C)=\frac{ \ker\left(\de:\ec^{p,q}(X)\to\ec^{p+1,q}(X)\right) \cap \ker\left(\dbar:\ec^{p,q}(X)\to\ec^{p,q+1}(X)\right) }{ \im\left(\ddbar:\ec^{p-1,q-1}(X)\to\ec^{p,q}(X)\right) }$$
\end{defn}

Il est facile de voir que $\bigoplus_{p,q}H^{p,q}_{BC}(X,\C)$ est muni d'une structure d'algèbre bigraduée induite par le produit extérieur des formes.\\

La définition même de ces différents espaces fournit des applications canoniques
$$H^{p,q}_{BC}(X,\C)\to H^{p+q}_{DR}(X,\C),\qquad H^{p,q}_{BC}(X,\C)\to H^{p,q}_{\dbar}(X,\C).$$

\subsection{Cas d'une variété kählerienne compacte}$\;$

Si $X$ est une variété kählerienne compacte, la théorie de Hodge fournit le résultat très utile suivant :
\begin{lemma}[du $\ddbar$] Soit $X$ une variété kählerienne compacte et $\alpha$ une forme de type $(p,q)$ $d$-fermée. Alors les conditions suivantes sont équivalentes :\\
$i)$ La forme $\alpha$ est $d$-fermée.\\
$ii)$ La forme $\alpha$ est $\de$-fermée.\\
$\overline{ii})$ La forme $\alpha$ est $\dbar$-fermée.\\
$iii)$ La forme $\alpha$ est $\ddbar$ fermée.\end{lemma}

On déduit de ce résultat que les applications suivantes sont des isomorphismes :
$$H^{p,q}_{BC}(X,\C)\stackrel{\sim}{\too}H^{p,q}_{\dbar}(X,\C),\quad \bigoplus_{p+q=k}H^{p,q}_{BC}(X,\C)\stackrel{\sim}{\too}H^k_{DR}(X,\C).$$
En particulier, la décomposition de Hodge est naturelle au sens où elle ne dépend pas du choix d'une métrique kählerienne sur $X$, de même que la symétrie de Hodge, puisque par définition on a bien sûr $H^{q,p}_{BC}(X,\C)=\overline{H^{p,q}_{BC}(X,\C)}$. 

\subsection{Exemple dans le cas non kählerien}$\;$

La situation est différente dans le cas non kählerien, comme le montre l'exemple de la variété d'Iwasawa. Soit $G$ le groupe de Heisenberg,  constitué des matrices de la forme
$$M(x,y,z)=\left(\begin{array}{ccc}1&x&z\\0&1&y\\0&0&1\end{array}\right),\quad x,y,z\in\C,$$
et soit $\Gamma$ le sous-groupe discret des matrices de la même forme dont les entrées sont des entiers de Gauss. On définit la variété d'Iwasawa $X$ comme le quotient $G/\Gamma$. Les formes holomorphes  $dx$, $dy$ et $dz-xdy$ sur $G$, obtenues comme composantes de $M^{-1}dM$, sont invariantes sous l'action à gauche de $\Gamma$. Elles induisent donc des formes holomorphes $\alpha,\beta,\gamma$ sur $X$. Le fait que la forme $\gamma$ n'est pas $d$-fermée (on a $\de\gamma=-\alpha\wedge\beta$) entraîne que $X$ n'est pas kählerienne. On peut facilement calculer les groupes de Dolbeault, De Rham et Bott-Chern à l'aide des formes $\alpha,\beta,\gamma$ (voir aussi dans [FG86] pour les deux premiers). 

\begin{prop} Le diamant de Hodge, les nombres de Betti et le diamant de Bott-Chern de la variété d'Iwasawa sont :
$$\begin{array}{ccccccccccccccc}
&&&1&&&&\qquad1\qquad&&&&1&&&\\
&&2&&3&&&4&&&3&&3&&\\
&2&&6&&3&&8&&2&&8&&2&\\
1&&6&&6&&1&10&1&&6&&6&&1\\
&3&&6&&2&&8&&3&&4&&3&\\
&&3&&2&&&4&&&2&&2&&\\
&&&1&&&&1&&&&1&&&\end{array}$$
\end{prop}

\begin{proof}
Les formes $\alpha,\beta,\gamma$ et leurs conjugués engendrent la $C^{\infty}(X)$-algèbre des formes différentielles sur $X$. Pour $M\in G$ on note $\tau_M$ l'opérateur de multiplication par $M$, bien défini sur $X$ et tel que $\tau_M^*\alpha=\alpha$, de même que pour $\beta$ et $\gamma$. Soit $d\mu$ une mesure invariante sur $X$ de volume 1, alors pour une forme $\omega$ sur $X$ la forme
$$\tilde{\omega}=\int_M\tau_M^*\omega\; d\mu(M)$$
est à coefficients constants en $\alpha,\beta,\gamma$ et leurs conjugués. De plus si $\omega$ est $d$-fermée (resp. $\dbar$-fermée) alors par homotopie $\tilde{\omega}$ lui est $d$-cohomologue (resp. $\dbar$-cohomologue). En conséquen\-ce, les groupes de De Rham et Dolbeault se calculent en testant la fermeture et l'exactitude des formes composées à partir de $\alpha,\beta,\gamma$. En particulier, on a une décomposition (qui n'est pas de Hodge) des groupes $H^{k}_{DR}(X,\C)$ selon le bidegré de ces générateurs. Par exemple,
$$\begin{array}{rl}
H^2_{DR}(X,\C)=&\bra \{\alpha\wedge\gamma\},\{\beta\wedge\gamma\} \ket_{\C}\\ 
&\oplus\bra\{\alpha\wedge\overline{\alpha}\},\{\alpha\wedge\overline{\beta}\},\{\beta\wedge\overline{\alpha}\},\{\beta\wedge\overline{\beta}\}  \ket_{\C}\\
&\oplus(\bra \{\overline{\alpha}\wedge\overline{\gamma}\},  \{\overline{\beta}\wedge\overline{\gamma}\}\ket_{\C}
\end{array}$$ 

Pour calculer les groupes de Bott-Chern, il faut tenir compte d'une part des formes $d$-fermées non $d$-exactes, données par les groupes de De Rham, et d'autre part des formes $d$-exactes mais non $\ddbar$-exactes (telles que $\alpha\wedge\beta$). Pour l'exemple, détaillons le calcul de $H^{2,2}_{BC}(X,\C)$ :\\

D'après le calcul de $H^4_{DR}(X,\C)$, une $(2,2)$-forme $d$-fermée est $d$-cohomologue à une combinaison linéaire des formes suivantes :
$$\alpha\wedge\gamma\wedge\overline{\alpha}\wedge\overline{\gamma},\; \alpha\wedge\gamma\wedge\overline{\beta}\wedge\overline{\gamma},\;\beta\wedge\gamma\wedge\overline{\alpha}\wedge\overline{\gamma},\;\beta\wedge\gamma\wedge\overline{\beta}\wedge\overline{\gamma}.$$
Il suffit donc de traiter le cas d'une forme $d$-exacte, que l'on notera $\omega=d\eta$. En décomposant $\eta$ selon son bidegré on obtient les équations :
$$\dbar\eta^{0,3}=0,\quad \de\eta^{0,3}+\dbar\eta^{1,2}=0,\quad\de\eta^{1,2}+\dbar\eta^{2,1}=\omega,\quad\de\eta^{2,1}+\dbar\eta^{3,0}=0,\quad\de\eta^{3,0}=0.$$

Si $\eta^{0,3}$ est $\dbar$-exact, disons $\eta^{0,3}=\dbar\xi$, alors $\omega=d(\eta-d\xi)$ où $\eta-d\xi$ n'a pas de composante de type $(0,3)$, donc on peut supposer $\eta^{0,3}=0$. On peut ainsi se ramener au cas où $\eta^{0,3}$ est une forme à coefficients constants en $\alpha,\beta,\gamma$, dont la classe est un générateur de $H^{0,3}_{\dbar}(X,\C)$. En l'occurence, ce groupe est de dimension 1, engendré par la classe de $\overline{\alpha}\wedge\overline{\beta}\wedge\overline{\gamma}$ qui est une forme $\de$-fermée. 
On déduit alors de la seconde équation que $\dbar\eta^{1,2}=0$. Si $\eta^{1,2}$ est $\dbar$-exact, le terme $\de\eta^{1,2}$ dans $\omega$ est $\ddbar$-exact et ne donne donc pas de contribution au groupe de Bott-Chern. On peut ainsi supposer que $\eta^{1,2}$ est une des formes à coefficients constants en $\alpha,\beta,\gamma$ générateur de $H^{1,2}_{\dbar}(X,\C)$, c'est-à-dire l'une des six formes suivantes :
$$\alpha\wedge\overline{\alpha}\wedge\overline{\gamma},\; \alpha\wedge\overline{\beta}\wedge\overline{\gamma},\;\beta\wedge\overline{\alpha}\wedge\overline{\gamma},\;\beta\wedge\overline{\beta}\wedge\overline{\gamma},\;\gamma\wedge\overline{\alpha}\wedge\overline{\gamma},\;\gamma\wedge\overline{\beta}\wedge\overline{\gamma}$$
Parmi celles-ci, les quatre premières sont $\de$-fermées mais pas les deux dernières ; on a en effet 
$$\de(\gamma\wedge\overline{\alpha}\wedge\overline{\gamma})=-\alpha\wedge\beta\wedge\overline{\alpha}\wedge\overline{\gamma}\quad\mbox{et}\quad\de(\gamma\wedge\overline{\beta}\wedge\overline{\gamma})=-\alpha\wedge\beta\wedge\overline{\beta}\wedge\overline{\gamma}$$
Ces deux éléments vont donner deux nouvelles classes dans $H^{2,2}_{BC}$. Par symétrie, la même étude concernant $\dbar\eta^{2,1}$ fournit les deux classes conjuguées, et on obtient finalement
$$\begin{array}{rcl}
H^{2,2}_{BC}(X,\C)=\;\bra\right.& \{\alpha\wedge\gamma\wedge\overline{\alpha}\wedge\overline{\gamma}\},\;\{\alpha\wedge\gamma\wedge\overline{\beta}\wedge\overline{\gamma}\},\;\{\beta\wedge\gamma\wedge\overline{\alpha}\wedge\overline{\gamma}\},\;\{\beta\wedge\gamma\wedge\overline{\beta}\wedge\overline{\gamma}\},&\\
&\{\alpha\wedge\beta\wedge\overline{\alpha}\wedge\overline{\gamma}\},\;\{\alpha\wedge\beta\wedge\overline{\beta}\wedge\overline{\gamma}\},\; \{\alpha\wedge\gamma\wedge\overline{\alpha}\wedge\overline{\beta}\},\;\{\beta\wedge\gamma\wedge\overline{\alpha}\wedge\overline{\beta}\}&\left.\ket_{\C},\end{array}$$
qui est de dimension $8$.
\end{proof}

Dans les tableaux suivants, on récapitule les générateurs des différents groupes de cohomologie considérés. On constate que pour cette variété l'application $H^{p,q}_{BC}(X,\C)\to H^{p,q}_{\dbar}(X,\C)$ est non surjective pour $p=1$ (par exemple, $\gamma$ n'a pas d'antécédent) et non injective pour $q=2$ (par exemple $\overline{\alpha}\wedge\overline{\beta}=\dbar(-\overline{\gamma})$ n'est pas $\ddbar$ exacte).

\begin{table}
$$\begin{array}{r||c|c|c|c|}\hline

& & \alpha\wedge\alphab\wedge\betab\wedge\gammab  & \alpha\wedge\gamma\wedge\alphab\wedge\betab\wedge\gammab  & \\
q=3& \alphab\wedge\betab\wedge\gammab  & \beta\wedge\alphab\wedge\betab\wedge\gammab  & \beta\wedge\gamma\wedge\alphab\wedge\betab\wedge\gammab  & \alpha\wedge\beta\wedge\gamma\wedge\alphab\wedge\betab\wedge\gammab  \\
&  & \gamma\wedge\alphab\wedge\betab\wedge\gammab  & \alpha\wedge\beta\wedge\alphab\wedge\betab\wedge\gammab  & \\\hline

& \alphab\wedge\gammab & \alpha\wedge\alphab\wedge\gammab\quad \alpha\wedge\betab\wedge\gammab & \alpha\wedge\gamma\wedge\alphab\wedge\gammab\quad \alpha\wedge\gamma\wedge\betab\wedge\gammab & \alpha\wedge\beta\wedge\gamma\wedge\alphab\wedge\gammab \\
q=2& \betab\wedge\gammab  & \beta\wedge\alphab\wedge\gammab\quad \beta\wedge\betab\wedge\gammab  &  \beta\wedge\gamma\wedge\alphab\wedge\gammab\quad \beta\wedge\gamma\wedge\betab\wedge\gammab  & \alpha\wedge\beta\wedge\gamma\wedge\betab\wedge\gammab  \\
&  &\gamma\wedge\alphab\wedge\gammab\quad \gamma\wedge\betab\wedge\gammab &  \alpha\wedge\beta\wedge\alphab\wedge\gammab\quad \alpha\wedge\beta\wedge\betab\wedge\gammab & \\\hline

& \alphab & \alpha\wedge\alphab\quad \alpha\wedge\betab & \alpha\wedge\gamma\wedge\alphab\quad \alpha\wedge\gamma\wedge\betab & \alpha\wedge\beta\wedge\gamma\wedge\alphab \\
q=1& \betab  & \beta\wedge\alphab\quad \beta\wedge\betab  &\beta\wedge\gamma\wedge\alphab\quad \beta\wedge\gamma\wedge\betab   &  \alpha\wedge\beta\wedge\gamma\wedge\betab  \\
& &\gamma\wedge\alphab\quad \gamma\wedge\betab & \alpha\wedge\beta\wedge\alphab\quad \alpha\wedge\beta\wedge\betab & \\\hline

& & \alpha & \alpha\wedge\gamma &\\
q=0& 1  & \beta  & \beta\wedge\gamma  & \alpha\wedge\beta\wedge\gamma \\
& & \gamma & \alpha\wedge\beta &\\\hline

&p=0&p=1&p=2&p=3\\\end{array}$$
\caption{Cohomologie de Dolbeault de la variété d'Iwasawa}
\end{table}

\begin{table}
$$\begin{array}{r||c|c|c|c|}\hline

q=3 & \alphab\wedge\betab\wedge\gammab & \alpha\wedge\alphab\wedge\betab\wedge\gammab  & \alpha\wedge\gamma\wedge\alphab\wedge\betab\wedge\gammab  & \alpha\wedge\beta\wedge\gamma\wedge\alphab\wedge\betab\wedge\gammab \\
&  & \alpha\wedge\alphab\wedge\betab\wedge\gammab  & \beta\wedge\gamma\wedge\alphab\wedge\betab\wedge\gammab  & \\\hline

q=2& \alphab\wedge\gammab & \alpha\wedge\alphab\wedge\gammab\quad \alpha\wedge\betab\wedge\gammab & \alpha\wedge\gamma\wedge\alphab\wedge\gammab\quad \alpha\wedge\gamma\wedge\betab\wedge\gammab & \alpha\wedge\beta\wedge\gamma\wedge\alphab\wedge\gammab \\
& \betab\wedge\gammab  & \beta\wedge\alphab\wedge\gammab\quad \beta\wedge\betab\wedge\gammab  &  \beta\wedge\gamma\wedge\alphab\wedge\gammab\quad \beta\wedge\gamma\wedge\betab\wedge\gammab  & \alpha\wedge\beta\wedge\gamma\wedge\betab\wedge\gammab  \\\hline

q=1& \alphab & \alpha\wedge\alphab\quad \alpha\wedge\betab & \alpha\wedge\gamma\wedge\alphab\quad \alpha\wedge\gamma\wedge\betab & \alpha\wedge\beta\wedge\gamma\wedge\alphab \\
& \betab  & \beta\wedge\alphab\quad \beta\wedge\betab  &\beta\wedge\gamma\wedge\alphab\quad \beta\wedge\gamma\wedge\betab   &  \alpha\wedge\beta\wedge\gamma\wedge\betab  \\\hline

q=0&1 & \alpha & \alpha\wedge\gamma &\alpha\wedge\beta\wedge\gamma  \\
&  & \beta  & \beta\wedge\gamma  & \\\hline

&p=0&p=1&p=2&p=3\\\end{array}$$
\caption{Cohomologie de De Rham de la variété d'Iwasawa}
\end{table}

\begin{table}
$$\begin{array}{r||c|c|c|c|}\hline

& & \alpha\wedge\alphab\wedge\betab\wedge\gammab  & \alpha\wedge\gamma\wedge\alphab\wedge\betab\wedge\gammab  & \\
q=3& \alphab\wedge\betab\wedge\gammab  & \beta\wedge\alphab\wedge\betab\wedge\gammab  & \beta\wedge\gamma\wedge\alphab\wedge\betab\wedge\gammab  & \alpha\wedge\beta\wedge\gamma\wedge\alphab\wedge\betab\wedge\gammab  \\
&  &  & \alpha\wedge\beta\wedge\alphab\wedge\betab\wedge\gammab  & \\\hline

& \alphab\wedge\gammab & \alpha\wedge\alphab\wedge\gammab\quad \alpha\wedge\betab\wedge\gammab & \alpha\wedge\gamma\wedge\alphab\wedge\gammab\quad \alpha\wedge\gamma\wedge\betab\wedge\gammab & \alpha\wedge\beta\wedge\gamma\wedge\alphab\wedge\gammab \\
q=2& \betab\wedge\gammab  & \beta\wedge\alphab\wedge\gammab\quad \beta\wedge\betab\wedge\gammab  &  \beta\wedge\gamma\wedge\alphab\wedge\gammab\quad \beta\wedge\gamma\wedge\betab\wedge\gammab  & \alpha\wedge\beta\wedge\gamma\wedge\betab\wedge\gammab  \\
&\alphab\wedge\betab  &\alpha\wedge\alphab\wedge\betab\quad \beta\wedge\alphab\wedge\betab &  \alpha\wedge\beta\wedge\alphab\wedge\gammab\quad \alpha\wedge\beta\wedge\betab\wedge\gammab & \alpha\wedge\beta\wedge\gamma\wedge\alphab\wedge\betab \\
&&&\alpha\wedge\gamma\wedge\alphab\wedge\betab\quad \beta\wedge\gamma\wedge\alphab\wedge\betab  &\\\hline

& \alphab & \alpha\wedge\alphab\quad \alpha\wedge\betab & \alpha\wedge\gamma\wedge\alphab\quad \alpha\wedge\gamma\wedge\betab & \alpha\wedge\beta\wedge\gamma\wedge\alphab \\
q=1& \betab  & \beta\wedge\alphab\quad \beta\wedge\betab  &\beta\wedge\gamma\wedge\alphab\quad \beta\wedge\gamma\wedge\betab   &  \alpha\wedge\beta\wedge\gamma\wedge\betab  \\
& & & \alpha\wedge\beta\wedge\alphab\quad \alpha\wedge\beta\wedge\betab & \\\hline

& & \alpha & \alpha\wedge\gamma &\\
q=0& 1  & \beta  & \beta\wedge\gamma  & \alpha\wedge\beta\wedge\gamma \\
& &  & \alpha\wedge\beta &\\\hline

&p=0&p=1&p=2&p=3\\\end{array}$$
\caption{Cohomologie de Bott-Chern de la variété d'Iwasawa}
\end{table}

\section{Théorie de Hodge de la cohomologie de Bott-Chern}

\subsection{Isomorphismes de Hodge classiques}$\;$

Soit $X$ une variété complexe compacte munie d'une métrique hermitienne $\omega$. Celle-ci permet de définir les adjoints des opérateurs $d$, $\partial$, $\dbar$ et les laplaciens associés :
$$\begin{array}{c}
\Delta=dd^*+d^*d:\ec^k(X)\to \ec^k(X),\\
\Delta'=\partial\partial^*+\partial^*\partial:\ec^{p,q}(X)\to \ec^{p,q}(X),\quad\Delta''=\dbar\,\dbar^*+\dbar^*\dbar:\ec^{p,q}(X)\to \ec^{p,q}(X).\end{array}$$
On note $\hc_{\Delta}^k(X)$ l'espace des formes harmoniques globales pour le laplacien $\Delta$ agissant sur les formes de degré $k$, etc. La théorie de Hodge fournit alors des isomorphismes
$$H^k_{DR}(X,\C)\cong\hc_{\Delta}^k(X),\quad H^{p,q}_{\dbar}(X,\C)\cong\hc_{\Delta''}^{p,q}(X).$$

Rappelons brièvement la démonstration du premier isomorphisme (celle du second est identique) : le laplacien $\Delta$ est un opérateur différentiel auto-adjoint, dont le symbole principal est donné, pour $x\in X$ et $\xi\in T^*_{X,x}$, par :
$$\sigma_{\Delta}(x,\xi)=-|\xi|^2_{\omega}\Id_{\Lambda^{p,q}T^*_{X,x}}.$$
C'est donc un opérateur elliptique, ce qui garantit alors une décomposition orthogonale
$$\ec^k(X)=\hc^k_{\Delta}(X)\oplus\im\Delta=\hc^k_{\Delta}(X)\oplus\im d\oplus\im d^*.$$
Etant donnée une $k$-forme $u$, elle s'écrit dans cette décomposition $u=h+dv+d^*w$. Sous cette écriture on a $du=0$ si et seulement si $dd^*w=0$, ce qui équivaut à $d^*w=0$. Finalement $\ker d=\hc^k_{\Delta}(X)\oplus\im d$, d'où l'isomorphisme de Hodge.

\subsection{Isomorphisme de Hodge pour la cohomologie de Bott-Chern}$\;$

Nous souhaitons obtenir un résultat analogue pour la cohomologie de Bott-Chern
$$H^{p,q}_{BC}(X,\C)=\dfrac{\ker \partial^{p,q}\cap \ker \dbar^{p,q}}{\ddbar A^{p-1,q-1}}.$$
Il semble donc naturel de considérer le ``laplacien de Bott-Chern'' suivant :
$$\Delta_{BC}^{p,q}=(\ddbar)(\ddbar)^*+\partial^*\partial+\dbar^*\dbar$$
En effet, c'est un opérateur auto-adjoint vérifiant, pour une forme $\fy$ de bidegré $(p,q)$,
$$\bra\bra\Delta_{BC}^{p,q}\fy,\fy\ket\ket=\|(\ddbar)^*\fy\|^2+\|\partial \fy\|^2+\|\dbar\fy\|^2.$$
Si la théorie précédente s'appliquait, on obtiendrait une décomposition orthogonale
$$\ec^{p,q}(X)=\hc^{p,q}_{\Delta_{BC}}(X)\oplus\im\ddbar\oplus(\im\partial^*+\im\dbar^*)$$
Pour une forme $u=h+\ddbar v+\partial^*w_1+\dbar^*w_2$, on aurait
$$\begin{array}{rcl}
\partial u=0\mbox{ et }\dbar u=0 & \eq & \partial(\partial^*w_1+\dbar^*w_2)=0\mbox{ et }\dbar(\partial^*w_1+\dbar^*w_2)=0\\
& \eq & \bra\!\bra \partial^*w_1+\dbar^*w_2, \partial^*w_1+\dbar^*w_2\ket\!\ket\\
& \eq & \partial^*w_1+\dbar^*w_2=0\end{array}$$
d'où l'on déduirait $\ker\partial^{p,q}\cap\ker\dbar^{p,q}=\hc_{\Delta_{BC}}^{p,q}(X)\oplus \im \ddbar$ et l'isomorphisme de Hodge souhaité.\\

Malheureusement :
\begin{prop} L'opérateur $\Delta_{BC}^{p,q}$ défini ci-dessus n'est pas elliptique.\end{prop}

\begin{proof} L'unique contribution au symbole principal est donné par le terme $\ddbar\,\dbar^*\partial^*$. Pour $x\in X$ et $\xi\in T^*_{X,x}$, exprimons les symboles principaux des différents opérateurs dans une base $(dz_I\wedge d\overline{z_J})_{|I|=p,|J|=q}$ associée à des coordonnées $(z_k)$ au voisinage de $x$.

\begin{notation}
Soit $I$ un multi-indice. Pour $i\in I$ on note $\alpha_i$ la position de $i$ dans $I$ et $I\moins i$ le multi-indice $I$ privé de $i$. pour $k\notin I$ on note $Ik$ le multi-indice $I\cup k$ réordonné et $\gamma_k$ la position de $k$ dans ce nouveau multi-indice.\\
Pour des entiers $i\not= k$ on note $\eps_{ik}$ le signe de $k-i$.\end{notation}

Les coefficients non nuls des matrices des symboles principaux considérés sont :

$$\begin{array}{cc}
(\sigma_{\partial}(x,\xi))_{I,J;Ik,J}=(-1)^{\gamma_k-1}\xi_k & (\sigma_{\partial^*}(x,\xi))_{I,J;I\moins i,J}=(-1)^{\alpha_i}\overline{\xi_k}\\
(\sigma_{\dbar}(x,\xi))_{I,J;I,Jl}=(-1)^{|I|+\delta_l-1}\overline{\xi_l} & (\sigma_{\dbar^*}(x,\xi))_{I,J;I,J\moins j}=(-1)^{|I|+\alpha_j}\xi_j\end{array}$$

On en déduit
$$(\sigma_{\ddbar\dbar^*\de^*}(x,\xi))_{I,J;I,J}=\left(\sum_{i\in I}|\xi_i|^2\right)\left(\sum_{j\in J}|\xi_j|^2\right),$$
$$(\sigma_{\ddbar\dbar^*\de^*}(x,\xi))_{I,J;Ik\moins i,J}=(-1)^{\alpha_i+\gamma_k-1}\eps_{ik}\overline{\xi_i}\xi_k\left(\sum_{j\in J}|\xi_j|^2\right),$$
$$(\sigma_{\ddbar\dbar^*\de^*}(x,\xi))_{I,J;I,Jl\moins j}=(-1)^{\beta_j+\delta_l-1}\eps_{jl}\xi_j\overline{\xi_l}\left(\sum_{i\in I}|\xi_i|^2\right),$$
$$(\sigma_{\ddbar\dbar^*\de^*}(x,\xi))_{I,J;Ik\moins i,Jl\moins j}=(-1)^{\alpha_i+\beta_j+\gamma_k+\delta_l}\eps_{ik}\eps_{jl}\overline{\xi_i}\xi_k\xi_k\overline{\xi_l}.$$

Ce symbole n'est pas injectif dès que $|I|\not= n$ ou $|J|\not=n$. S'il existe par exemple $k\notin I$ alors $\sigma_{\ddbar\dbar^*\de^*}(x,dz_k)=0$.
\end{proof}

On va donc modifier l'opérateur $\Delta_{BC}^{p,q}$ en ajoutant d'autres termes d'ordre 4 pour le rendre elliptique. Concrètement, nous posons
$$\widetilde{\Delta}_{BC}^{p,q}=\ddbar\,\dbar^*\partial^*+\dbar^*\partial^*\partial\dbar+\dbar^*\partial\partial^*\dbar+\partial^*\dbar\,\dbar^*\partial+\dbar^*\dbar+\partial^*\partial$$
Dans le calcul de la matrice du symbole principal associé, tous les termes non diagonaux se compensent. Quant aux termes diagonaux, on a :

$$(\sigma_{\dbar^*\de^*\ddbar}(x,\xi))_{I,J;I,J}=\left(\sum_{k\notin I}|\xi_k|^2\right)\left(\sum_{l\notin J}|\xi_l|^2\right)$$
$$(\sigma_{\dbar^*\de\de^*\dbar}(x,\xi))_{I,J;I,J}=\left(\sum_{i\in I}|\xi_i|^2\right)\left(\sum_{l\notin J}|\xi_l|^2\right)$$
$$(\sigma_{\de^*\dbar\dbar^*\de}(x,\xi))_{I,J;I,J}=\left(\sum_{k\notin I}|\xi_k|^2\right)\left(\sum_{j\in J}|\xi_j|^2\right)$$

Ainsi,
$$\sigma_{\widetilde{\Delta}_{BC}^{p,q}}(x,\xi)=|\xi|_{\omega}^4\Id_{\Lambda^{p,q}T^*_{X,x}}$$
L'opérateur $\widetilde{\Delta}_{BC}^{p,q}$ est donc bien elliptique, et par ailleurs a le même noyau que $\Delta_{BC}^{p,q}$ :
$$\begin{array}{rcl}
\widetilde{\Delta}_{BC}^{p,q}u=0 & \eq & \dbar u=\partial u=\dbar^*\partial^*u=\ddbar u=\partial^*\partial u=\dbar^*\partial u=0\\
& \eq & \dbar u=\partial u=\dbar^*\partial^*u=0.\end{array}$$

On a donc démontré :
\begin{thm}[Décomposition et isomorphisme de Hodge pour la cohomologie de Bott-Chern]
On a une décomposition orthogonale
$$\ec^{p,q}(X)=\hc^{p,q}_{\Delta_{BC}}(X)\oplus\im\ddbar\oplus(\im\partial^*+\im\dbar^*)$$
et un isomorphisme
$$H^{p,q}_{BC}(X,\C)\cong\hc^{p,q}_{\Delta_{BC}}(X).$$
\end{thm}

\begin{cor} Si $X$ est une variété complexe compacte alors ses groupes de cohomologie de Bott-Chern sont de dimension finie.
\end{cor}
\begin{proof} Le noyau d'un opérateur différentiel elliptique sur une variété compacte est de dimension finie.
\end{proof}

Terminons par le cas kählerien : l'utilisation d'identités de commutations telles que $[\dbar^*,\omega\wedge\cdot]=i\dbar$ donne le résultat classique suivant, et son analogue de Bott-Chern :
\begin{prop} Si $(X,\omega)$ est kählerienne compacte, on a les identités
$$\Delta'=\Delta''=\frac{1}{2}\Delta,\quad \widetilde{\Delta}_{BC}=\Delta''\Delta''+\partial^*\partial+\dbar^*\dbar$$
et les espaces de formes harmoniques $\hc_{\Delta}^{p+q}\cap \ec^{p,q}(X)$, $\hc_{\Delta''}^{p,q}$ et $\hc_{\widetilde{\Delta}_{BC}}^{p,q}$ coïncident.\\ 
De plus, la forme de Kähler $\omega\in \ec^{1,1}(X)$ est harmonique pour tous les laplaciens considérés.
\end{prop}

\subsection{Cohomologie d'Aeppli}$\;$

L'exemple déjà considéré de la variété d'Iwasawa montre qu'il n'y a pas d'analogue à la dualité de Serre pour la cohomologie de Bott-Chern (cette dualité correspond à une symétrie centrale du diamant). Cela provient de la dissymétrie des opérateurs utilisés dans la définition des groupes de Bott-Chern.

Pour une variété compacte $X$, le dual de $H^{p,q}_{BC}(X,\C)$ est le groupe de cohomologie d'Aeppli $H^{n-p,n-q}_A(X,\C)$ où
$$H^{p,q}_{A}(X,\C)=\frac{\ker\ddbar:\ec^{p,q}(X)\to\ec^{p+1,q+1}(X)}{(\im \partial:\ec^{p-1,q}(X)\to\ec^{p,q}(X))+(\im \dbar:\ec^{p,q-1}(X)\to\ec^{p,q}(X))}.$$

\begin{lemma} Le produit extérieur induit une application bilinéaire
$$\wedge : H^{p,q}_{BC}(X,\C)\times H^{r,s}_A(X,\C)\to H^{p+r,q+s}_A(X,\C).$$
\end{lemma}
\begin{proof} Le produit d'une forme $d$-fermée et d'une forme $\ddbar$-fermée est une forme $\ddbar$-fermée. Par ailleurs, le produit d'une forme $d$-fermée par une forme $d$-exacte est $d$-exact, et le produit d'une forme $\ddbar$-exacte par une forme $\ddbar$-fermée est lui aussi $d$-exact. Pour montrer ce dernier point, on utilise la formule suivante : si $\alpha$ est une $(p-1,q-1)$-forme et $\beta$ une $(r,s)$-forme $\ddbar$-fermée alors
$$\ddbar\alpha\wedge\beta=\frac{1}{2}d\left[(\dbar\alpha-\partial\alpha)\wedge\beta + (-1)^{p+q}\alpha\wedge (\partial\beta-\dbar\beta)\right].$$\end{proof}

En particulier on a un appariement
$$H^{p,q}_{BC}(X,\C)\times H^{n-p,n-q}_A(X,\C)\to H^{n,n}_A(X,\C)=H^{n,n}_{DR}(X,\C)\stackrel{\sim}{\to}\C,$$
la dernière flèche étant donnée par l'intégration sur $X$.\\

De la même manière que les groupes de Bott-Chern, les groupes de cohomologie d'Aeppli d'une variété compacte s'identifient aux formes harmoniques pour un laplacien d'ordre 4. L'opérateur a priori naturel
$$\Delta_A^{p,q}=\partial\partial^*+\dbar\,\dbar^*+(\ddbar)^*(\ddbar)$$
n'étant pas elliptique, on le remplace par l'opérateur elliptique de même noyau
$$\widetilde{\Delta}_A^{p,q}=\partial\partial^*+\dbar\,\dbar^*+\dbar^*\partial^*\ddbar+\ddbar\dbar^*\partial^*+\partial\dbar^*\dbar\partial^*+\dbar\partial^*\partial\dbar^*.$$
La théorie des opérateurs elliptiques fournit une décomposition orthogonale
$$\ec^{p,q}(X)=\hc_{\tilde{\Delta}_A}^{p,q}(X)\oplus (\im\partial+\im\dbar)\oplus\im(\ddbar)^*$$
d'où l'on déduit l'isomorphisme de Hodge $\displaystyle{H^{p,q}_A(X,\C)\cong \hc_{\tilde{\Delta}_A}^{p,q}(X)}$ de la même manière que nous l'avons fait pour les groupes de Bott-Chern.\\

Les isomorphismes de Hodge ainsi construits, outre qu'ils prouvent la finitude des groupes de Bott-Chern et d'Aeppli, permettent également de démontrer la dualité entre ces espaces. En effet, en notant $*:\ec^{p,q}(X)\to \ec^{n-p,n-q}(X)$ l'opérateur de Hodge défini par $u\wedge *v=\bra u,v\ket dV$, on a
$$\partial^*=\pm *\partial *,\quad \dbar^*=\pm *\dbar *,$$ 
d'où l'on déduit que
$$\begin{array}{rcl}
u\in\hc^{p,q}_{\widetilde{\Delta}_{BC}}(X) & \eq & \partial u=0,\;\dbar u=0,\;(\ddbar)^*u=0\\
& \eq & \dbar^*(* u)=0,\;\partial^*(*u)=0,\;\ddbar(*u)=0\\
& \eq & *u\in \hc^{n-p,n-q}_{\widetilde{\Delta}_A}(X)\end{array}$$
L'opérateur de Hodge réalise donc un isomorphisme entre $H^{p,q}_{BC}(X,\C)$ et $H^{n-p,n-q}_A(X,\C)$.

\begin{rem}En particulier, dans le cas compact, $H^{0,0}_A(X,\C)$ est l'espace des fonctions $\ddbar$-fermées, donc constantes, ce qui prouve que l'on a toujours $H^{n,n}_{BC}(X,\C)=\C$ pour une variété compacte.\end{rem}

\section{Applications à la théorie des déformations}

\subsection{Le théorème de Kodaira et Spencer}$\;$

Nous nous intéressons ici au célèbre théorème affirmant qu'une petite déformation d'une variété kählerienne est kählerienne. La preuve de ce théorème par Kodaira et Spencer, telle qu'on peut la trouver dans [MK], repose sur l'étude d'un opérateur elliptique d'ordre 4 qui n'est autre que l'opérateur ``laplacien de Bott-Chern'' précédemment défini. Ceci permet un éclairage de la preuve en termes de cohomologie de Bott-Chern.\\

Commençons par l'énoncé précis du théorème. Celui-ci s'applique  à une petite déforma\-tion différentiable de variétés complexes compactes. On entend par là une submersion $C^{\infty}$ : $\pi:\xc\too S$ où $S$ est une boule dans $\R^m$ et $\xc$ une variété différentiable, telle que chaque fibre de $\pi$ soit une variété complexe compacte. 
Plus précisément, $\xc$ est recouvert par des ouverts de carte de coordonnées $(z_1,\ldots,z_n,t_1,\ldots,t_m)\in U\subset \C^n \times\R^m$ telles que pour $t=(t_1,\ldots,t_m)$ fixé, $(z_1,\ldots,z_n)$ soient des coordonnées holomorphes sur $X_t=\pi^{-1}(t)$. Sous ces conditions, toutes les variétés $X_t$ sont $C^{\infty}$-difféomorphes entre elles.

\begin{thm}[Kodaira, Spencer]
Soit  $\pi:\xc\too S$ une déformation différentiable de variétés complexes compactes. On suppose que la fibre centrale $X_0$ est kählerienne. Alors il existe $\eps>0$ tel que pour $\|t\|<\eps$, la variété $X_t$ est kählerienne.
\end{thm}

\begin{proof}
Soit $\omega_0$ une métrique de Kähler sur $X_0$. A l'aide d'une partition de l'unité, on l'étend en une $2$-forme différentielle $\Omega$ sur $\xc$ telle que $\omega_t:=\Omega|_{X_t}$ soit une métrique hermitienne sur $X_t$.

Cette métrique hermitienne permet de construire des adjoints et par suite des opérateurs $\Delta''_t$ et $\widetilde{\Delta}_{BC,t}$ agissant sur les formes sur $X_t$. De plus, la métrique $\omega_t$ variant de manière différentiable, ces opérateurs varient de manière différentiable avec $t$. 

Par des résultats classiques de théorie spectrale, on en déduit que les valeurs propres de ces opérateurs, ordonnées en tenant compte des multiplicités, varient continument avec $t$. En particulier, pour $t$ suffisamment proche d'un $t_0$ fixé, le noyau de $\Delta''_t$ (resp. $\widetilde{\Delta}_{BC,t}$) ne peut avoir une dimension supérieure à celle du noyau de $\Delta''_{t_0}$ (resp. $\widetilde{\Delta}_{BC,t_0}$). Autrement dit, la dimension des espaces de formes harmoniques est une fonction semi-continue supérieurement de $t$. 
En combinant ce résultat avec les isomorphismes de Hodge, on obtient le premier ingrédient de la preuve :

\begin{lemma} Les nombres de Hodge $h^{p,q}_{\dbar}(X_t)$ et de Bott-Chern $h^{p,q}_{BC}(X_t)$ sont des fonctions semi-continues supérieurement de $t$.\end{lemma}

Remarquons que les nombres de Betti $b^k(X)$, eux, sont constants puisqu'il ne dépendent que de la structure différentielle de la variété. 

Le deuxième ingrédient est une inégalité liant certains nombres de Betti, Hodge et Bott-Chern :
\begin{lemma} Pour toute variété complexe compacte on a l'inégalité
$$b^2(X)\les 2 h^{0,2}_{\dbar}(X)+h^{1,1}_{BC}(X).$$
Si $X$ est kählerienne compacte, il y a égalité.
\end{lemma}
\begin{proof} On considère la suite exacte suivante :
$$0 \to \dfrac{\ker d^{1,1}\cap dA^1}{\ddbar A^0} \to H^{1,1}_{BC}(X,\C) \to  H^2_{DR}(X,\C) \stackrel{u}{\to} \overline{H^{0,2}_{\dbar}(X,\C)}\oplus H^{0,2}_{\dbar}(X,\C) \to \mathrm{coker(u)}\to 0$$
où l'application $u$ envoie une $2$-forme $\psi$ sur ses composantes de type $(2,0)$ et $(0,2)$.

On en déduit
$$h^{1,1}_{BC}(X)  = \dim \dfrac{\ker d^{1,1}\cap dA^1}{\ddbar A^0}+b^2(X)+\dim \mathrm{coker} u-2 h^{0,2}_{\dbar}(X)\ges b^2(X)-2h^{0,2}_{\dbar}(X).$$

Le fait que l'inégalité soit une égalité dans le cas kählerien résulte de l'égalité entre les groupes de Dolbeault et les groupes de Bott-Chern et de la décomposition de Hodge : si $X$ est kählerienne,
$$h^{0,2}_{\dbar}(X)+h^{1,1}_{BC}(X)=h^{2,0}_{\dbar}(X)+h^{0,2}_{\dbar}(X)+h^{1,1}_{\dbar}(X)=b^2(X).$$
\end{proof}

La variété $X_0$ étant kählerienne, l'inégalité
$$b^2\les 2h^{0,2}_{\dbar}(X_t)+h^{1,1}_{BC}(X_t)$$ 
est une égalité pour $t=0$. Le fait que $h^{0,2}_{\dbar}(X_t)$ et $h^{1,1}_{BC}(X_t)$ soient des fonctions semi-continue supérieurement de $t$ implique alors que ces deux nombres sont constants pour $t$ proche de $0$. 
Dans ce cas, l'espace $\hc^{1,1}_{\widetilde{\Delta}_{BC}}(X_t)$ varie différentiablement avec $t$, de même que la projection orthogonale $h_t$ de $A^{1,1}(X_t)$ sur cet espace.

En posant alors
$$\widetilde{\omega}_t=\frac{1}{2}(h_t\omega_t+\overline{h_t\omega_t}),$$
on obtient sur $X_t$ une $(1,1)$-forme réelle $d$-fermée, vérifiant
$$\lim_{t\to 0}\widetilde{\omega}_t=\frac{1}{2}(h_0\omega_0+\overline{h_0\omega_0})=\omega_0$$
car $h_0\omega_0=\omega_0$ puisque $\omega_0$ est harmonique.
Comme $\omega_0$ est définie-positive, il en va de même de $\widetilde{\omega}_t$ pour $t$ assez petit, et on a bien construit une métrique de Kähler sur $X_t$.
\end{proof}

\subsection{Déformations de la variété d'Iwasawa}

$\;$
\bigskip

\noindent
(Note de J.-P. Demailly, septembre 2007~: Michel Schweitzer, qui a
calculé explicitement ces déformations, n'a pas eu le temps de rédiger
cette partie après son départ de l'Université de Grenoble en juillet
2007 ; contacter directement l'auteur pour plus de détails).

\vfill

\pagebreak

\pagestyle{fancy}
\renewcommand{\headrulewidth}{0.4pt} 
\renewcommand{\footrulewidth}{0pt}
\lfoot{} \cfoot{} \rfoot{} 
\lhead[\thepage / Cohomologie de Bott-Chern entière]{\textit{Autour de la cohomologie de Bott-Chern}} \chead{} 
\rhead[\textit{Autour de la cohomologie de Bott-Chern}]{Cohomologie de Bott-Chern entière / \thepage}

\part*{Deuxième partie : Cohomologie de Bott-Chern entière}

\section{Interprétation hypercohomologique de la cohomologie de Bott-Chern}

Dans ce paragraphe on construit un isomorphisme entre les groupes de cohomologie de Bott-Chern et des groupes d'hypercohomologie de complexes définis à l'aide des faisceaux de formes holomorphes et antiholomorphes. Ce procédé, utilisé par exemple dans [Dem93] pour prouver la finitude de ces groupes, nous sert ici à définir la cohomologie de Bott-Chern entière.

\subsection{Lemme de résolubilité locale}
On rassemble dans la proposition suivante les liens entre la résolubilité locale des différents opérateurs $d,\de,\dbar,\ddbar$, conséquences des lemmes de Poincaré et de Dolbeault-Grothendieck :

\begin{lemma}[de résolubilité locale]On se place sur une boule $U\subset\C^n$.\\
1. Soit $\theta$ une forme de degré $k\ges 1$ de composantes de type $(p,q)$ nulles sauf si $p_1\les p\les p_2$ ($p_1<p_2$). Si $\theta$ est $d$-fermée alors $\theta=d\alpha$ avec $\alpha$ de degré $k-1$, de composantes de type $(p,q)$ nulles sauf si $p_1\les p\les p_2-1$.\\
2. Soit $\theta$ une $(p,q)$-forme, $d$-fermée : $\de\theta=\dbar\theta=0$.\\
$i)$ (Lemme du $\ddbar$ local) Si $p\ges 1$ et $q\ges 1$ alors $\theta\in\ddbar\ec^{p-1,q-1}(U)$;\\
$ii)$ (Lemme de Poicaré holomorphe) Si $p\ges 1$ et $q=0$ alors $\theta\in\de\Omega^{p-1}(U)$;\\
$\overline{ii})$ Si $p=0$ et $q\ges 1$ alors $\theta\in \dbar\ombar{q-1}(U)$;\\
$iii)$ Si $p=q=0$ alors $\theta$ est constante.\\
3. Soit $\theta$ une $(p,q)$-forme, $\ddbar$-fermée : $\ddbar\theta=0$. Alors $\theta$ est somme d'une forme $\de$-fermée et d'une forme $\dbar$-fermée. Autrement dit :\\
$i)$ Si $p\ges 1$ et $q\ges 1$ alors $\theta\in\dbar\ec^{p,q-1}(U)+\de\ec^{p-1,q}(U)$;\\
$ii)$ Si $p\ges 1$ et $q=0$ alors $\theta\in\Omega^{p}(U)+\de\ec^{p-1,0}(U)$;\\
$\overline{ii})$ Si $p=0$ et $q\ges 1$ alors $\theta\in \dbar\ec^{0,q-1}(U)+\ombar{q}(U)$;\\
$iii)$ Si $p=q=0$ alors $\theta\in\oc(U)+\ocbar(U)$.\\
4. Soit $\theta$ une forme de degré $k\ges 1$ que l'on suppose ``presque'' $d$-fermée, c'est-à-dire que l'on n'impose pas $\dbar\theta^{p_1,q_1}=0$ si $p_1\ges 1$ ni $\de\theta^{p_2,q_2}=0$ si $q_2\ges 1$.
Alors il existe des formes $\gamma^{p_1,q_1},\alpha^{p_1,q_1-1},\ldots,\alpha^{p_2-1,q_2},\gamma^{p_2,q_2}$ telles que $\gamma^{p_1,q_1}$ est $\de$-fermée, $\gamma^{p_2,q_2}$ est $\dbar$-fermée, et
$$\begin{array}{c}
\theta^{p_1,q_1}=\gamma^{p_1,q_1}+\dbar\alpha^{p_1,q_1-1},\;\theta^{p_1+1,q_1-1}=\de\alpha^{p_1,q_1-1}+\dbar\alpha^{p_1+1,q_1-2},\ldots,\\
\theta^{p_2-1,q_2+1}=\de\alpha^{p_2-2,q_2+1}+\dbar\alpha^{p_2-1,q_2},\;\theta^{p_2,q_2}=\de\alpha^{p_2-1,q_2}+\gamma^{p_2,q_2}.\end{array}$$
\end{lemma}

\begin{proof}
1. On écrit $\theta=d\beta$ (lemme de Poincaré). Si $p_1=0$ et $p_2=k$ il n'y a rien à prouver, on peut donc supposer $k\ges 2$ et $p_1>0$. On déduit alors de $\theta=d\beta$ que $\dbar\beta^{0,k-1}=0$. D'après le lemme de Dolbeault-Grothendieck, on peut écrire $\beta^{0,k-1}=\dbar\gamma^{0,k-2}$. On pose $\tilde{\beta}=\beta-d\gamma^{0,k-2}$, on a toujours $\theta=d\tilde{\beta}$, mais cette fois $\tilde{\beta}^{0,k-1}=0$. On peut donc supposer que $\beta$ n'a pas de composante de type $(0,k-1)$. En répétant ce raisonnement, on voit que l'on peut supposer que $\beta$ n'a pas de composante de type $(p,q)$ pour $p<p_1$. De même, si $p_2<k$, on élimine les composantes de type $(p,q)$ pour $p>p_2$ en utilisant le lemme de Dolbeault-Grothendieck anti-holomorphe.

2. Le $iii)$ est évident, et $ii)$ et $\overline{ii})$ sont équivalents, donc on peut supposer $p\ges 1$. On applique le $1.$ à la forme $\tilde{\theta}$ définie par $\tilde{\theta}^{p-1,q+1}=0$ et $\tilde{\theta}^{p,q}=\theta$ : il existe $\alpha$ de type pur $(p-1,q)$ telle que $\tilde{\theta}=d\alpha$, donc $\theta=\de\alpha$ avec $\dbar\alpha=0$, ce qui montre le résultat.

3. Posons $\theta^{p+1,q}=\de\theta^{p,q}$ : c'est une forme $d$-fermée, donc d'après le $2.$, $\theta^{p+1,q}=\de\alpha^{p,q}$ avec $\dbar\alpha^{p,q}=0$. Ainsi $\theta=(\theta-\alpha)+\alpha$ est somme d'une forme $\de$-fermée et d'une forme $\dbar$-fermée.

4. On applique le 3. à $\theta^{p_1,q_1}$ : on peut écrire $\theta^{p_1,q_1}=\gamma^{p_1,q_1}+\dbar\alpha^{p_1,q_1-1}$. Puis
$$\dbar(\theta^{p_1+1,q_1-1}-\de\alpha^{p_1,q_1-1})=\de(-\theta^{p_1,q_1}+\dbar\alpha^{p_1,q_1-1})=0$$
donc $\theta^{p_1+1,q_1-1}=\de\alpha^{p_1,q_1-1}+\gamma^{p_1+1,q_1-1}$ avec $\gamma^{p_1+1,q_1-1}$ $\dbar$-fermée, et on continue.
\end{proof}

\begin{rem} Puisque les lemmes de Poincaré et de Dolbeault-Grothendieck sont vrais pour les courants, on a les mêmes résultats en remplaçant les espaces de formes par les espaces de courants.\end{rem}

\subsection{Complexes associés aux groupes de Bott-Chern}
On fixe $p\ges 1,q\ges 1$, et on définit le complexe de faisceaux $\lc^{\bullet}_{p,q}$ (qui sera le plus souvent noté $\lc^{\bullet}$) par
$$\lc^k_{p,q}=\bigoplus_{\substack{r+s=k \\ r<p,s<q}}\ec^{r,s}\quad \mbox{si }k\les p+q-2,$$
$$\lc^{k-1}_{p-1,q-1}=\bigoplus_{\substack{r+s=k \\ r\ges p,s\ges q}}\ec^{r,s}\quad \mbox{si }k\ges p+q,$$
et la différentielle :
$$\lc^0\stackrel{\pr_{\lc^1}\circ d}{\too}\lc^1\stackrel{\pr_{\lc^2}\circ d}{\too}\ldots\to\lc^{k-2}\stackrel{\ddbar}{\too}\lc^{k-1}\stackrel{d}{\too}\lc^k\stackrel{d}{\too}\ldots$$

On a par construction 
$$H^{p,q}_{BC}(X,\C)=H^{p+q-1}(\lc^{\bullet}(X))\cong \Hb^{p+q-1}(X,\lc^{\bullet}),$$
l'isomorphisme venant du fait que les faisceaux $\lc^k$ sont mous.

On définit un sous-complexe $\sca^{\bullet}=\sca_{p,q}^{\bullet}$ de $\lc^{\bullet}=\lc_{p,q}^{\bullet}$ par :
$$({\sca'}_p^{\bullet},\de):\; \oc\to\Omega^1\to\ldots\to\Omega^{p-1}\to 0,\qquad ({\sca''}_q^{\bullet},\dbar):\; \ocbar\to\ombar{1}\to\ldots\to\ombar{q-1}\to 0$$
$$\sca_{p,q}^{\bullet}={\sca'}_p^{\bullet}+{\sca''}_q^{\bullet}:\; \oc+\ocbar\to\Omega^1\oplus\ombar{1}\to\ldots\Omega^{p-1}\oplus\ombar{p-1}\to\ombar{p}\to\ldots\to\ombar{q-1}\to 0$$

\begin{prop} L'inclusion $\sca^{\bullet}\subset\lc^{\bullet}$ est un quasi-isomorphisme.\end{prop}
\begin{proof} C'est la retraduction du lemme de résolubilité locale, qui calcule les faisceaux de cohomologie $\hc^k(\lc^{\bullet})$ et $\hc^k(\sca^{\bullet})$ : si $k\ges\max(p,q)$, une forme $d_{\lc}$-fermée est localement $d_{\lc}$-exacte (par la partie 1. si $k\ges p+q$, la partie 2. si $k=p+q-1$, la partie 3. si $k=p+q-2$ et la partie 4. si $k<p+q-2$).\\ 
Soit $k=p-1\ges q$, et soit $\theta\in\lc^k(U)$ une section $d_{\lc}$-fermée sur une boule $U\subset X$, de composantes $\theta^{p-q,q-1},\ldots,\theta^{p-1,0}$. D'après la partie 4. du lemme, il existe sur $U$ une $p-2$-forme $\alpha=\alpha^{p-q-1,q-1}+\ldots+\alpha^{p-2,0}$ et une $p-1$ forme holomorphe $u^{p-1,0}$ telles que $\theta=d_{\lc}\alpha+u^{p-1,0}$. Si l'on modifie $\theta$ par une forme $d_{\lc}$-exacte, $u^{p-1,0}$ est modifié par une forme holomorphe $\de$-fermée, d'où l'on déduit que
$$\hc^{p-1}(\lc^{\bullet})=\hc^{p-1}(\sca^{\bullet})=\Omega^{p-1}/\de\Omega^{p-2}.$$
Si $p-1>k\ges q-1$, le même raisonnement s'applique, à ceci près que la forme holomorphe $u^{k,0}$ est $\de$-fermée, et donc 
$$\hc^{k}(\lc^{\bullet})=\hc^{k}(\sca^{\bullet})=0.$$
Les autres cas se traitent de la même manière, excepté pour $k=0$ où il est clair que $\hc^0(\lc^{\bullet})=\hc^0(\sca^{\bullet})=\oc+\ocbar$. \end{proof}
\begin{rem} Pour une preuve formelle de ce résultat, utilisant des suites spectrales, on renvoie au livre de J.-P. Demailly [Dem93], paragraphe VI.12.\end{rem}

Le complexe $\sca^{\bullet}$ est peu pratique du fait de la somme non directe $\oc+\ocbar$. On le modifie donc légèrement en définissant le complexe $\bc_{p,q}=\bc^{\bullet}_{p,q}$ par :
$$\bc_{p,q}^{\bullet} :\; \C\stackrel{(+,-)}{\too}\oc\oplus\ocbar\to\Omega^1\oplus\ombar{1}\to\ldots\Omega^{p-1}\oplus\ombar{p-1}\to\ombar{p}\to\ldots\to\ombar{q-1}\to 0.$$
\begin{prop}
L'application naturelle de $\bc^{\bullet}$ dans $\sca^{\bullet}[1]$ :
$$\begin{array}{ccccccc}
\C&\stackrel{(+,-)}{\too}&\oc\oplus\ocbar&\to&\Omega^1\oplus\ombar{1}&\to&\ldots\\
\downarrow&&\downarrow +&&\downarrow&&\\
0&\too&\oc +\ocbar&\to&\Omega^1\oplus\ombar{1}&\to&\ldots
\end{array}$$
est un quasi-isomorphisme.\end{prop}
\begin{proof} La seule vérification à faire est en degré 1 :
$$\hc^1(\bc^{\bullet})=(\C\oplus\C)/\C(1,-1)\stackrel{+}{\too}\C=\hc^1(\sca^{\bullet})$$
est un isomorphisme.\end{proof}

Les trois complexes $\lc^{\bullet}[1]$, $\sca^{\bullet}[1]$ et $\bc^{\bullet}$ étant quasi-isomorphes, ils ont même hypercohomologie. On a en particulier la suite d'isomorphismes
$$H^{p,q}_{BC}(X,\C)\cong\Hb^{p+q}(X,\lc_{p,q}^{\bullet}[1])\cong\Hb^{p+q}(X,\sca_{p,q}^{\bullet}[1])\cong \Hb^{p+q}(X,\bc_{p,q}^{\bullet}),$$
et c'est cette dernière écriture que nous utiliserons par la suite.

\subsection{Cas où $p=0$ ou $q=0$}
Supposons par exemple que $p>0$ et $q=0$. On peut définir de la même manière le complexe $\lc_{p,0}^{\bullet}$, qui devient simplement
$$\lc_{p,0}^k=0\quad \mbox{si }k\les p-2,$$
$$\lc_{p,0}^{k-1}=\bigoplus_{\substack{r+s=k \\ r\ges p}}\ec^{p,q}\quad \mbox{si }k\ges p,$$
avec la différentielle usuelle, et on a toujours $H^{p,0}_{BC}(X,\C)=\Hb^{p_0-1}(X,\lc^{\bullet})$. Il n'y a par contre pas d'équivalent au complexe $\sca^{\bullet}$ : il est clair que
$$H^{p_0,0}_{BC}(X)\not=\Hb^{p_0-1}(X,\oc\to\ldots\to\Omega^{p_0-1})=\Hb^{p-1}(X,\Omega^{\bullet})=H^{p_0-1}(X,\C).$$
En revanche, le complexe $\bc_{p,0}^{\bullet}$ calcule bien la cohomologie de Bott-Chern :
\begin{prop} L'application naturelle de $\bc_{p,0}^{\bullet}$ dans $\lc_{p,0}^{\bullet}[1]$ :
$$\begin{array}{ccccccccccc}
\C&\to&\oc&\to&\ldots&\to&\Omega^{p-1}&\to& 0&&\\
\downarrow&&\downarrow&&&&\de\downarrow&&\downarrow&&\\
0&\to&0&\to&\ldots&\to&\ec^{p,0}&\to&\ec^{p,1}\oplus\ec^{p+1,0} &\to&\ldots
\end{array}$$
est un quasi-isomorphisme.
\end{prop}
\begin{proof}A nouveau, le lemme de résolubilité locale calcule les faisceaux $\hc^{k}(\lc^{\bullet}[1])$ et $\hc^{k}(\bc^{\bullet})$, qui sont nuls si $k\not=p$. Si $k=p\ges 2$, on obtient
$$\hc^{p}(\bc^{\bullet})=\Omega^{p-1}/\de\Omega^{p-2}\stackrel{\de}{\too}\de\Omega^{p-1}=\ker(\de:\Omega^{p}\to\ec^{p+1})=\hc^{p}(\lc^{\bullet}[1])$$
qui est un isomorphisme ; si $k=p=1$,
$$\hc^{1}(\bc^{\bullet})=\oc/\C\stackrel{\de}{\too}\de\oc=\ker(\de:\Omega^{1}\to\ec^{2})=\hc^{1}(\lc^{\bullet}[1]).$$
est également un isomorphisme.
\end{proof}

Le cas $p=0$, $q>0$ est symétrique, mais cette fois on choisit $-\dbar$ comme quasi-isomorphisme ; enfin, le cas $p=q=0$ rentre dans cette description puisque $H^{0,0}_{BC}(X,\C)=H^0(X,\C)$. 

\subsection{Conséquences de l'interprétation hypercohomologique}

Commençons par trois remarques importantes relatives à ce qui précède :

Les lemmes de Poincaré et de Dolbeault-Grothendieck étant vrais lorsque l'on remplace les formes $C^{\infty}$ par des courants, le lemme de résolubilité locale et l'interprétation hypercohomologique qui en découle restent vrais pour des courants.

Par ailleurs, cette interprétation hypercohomologique montre que pour une variété compacte, les groupes de cohomologie de Bott-Chern sont de dimension finie : il suffit pour cela de raisonner par récurrence sur le complexe $\sca$ en utilisant le fait que les groupes $H^q(X,\Omega^p)$ sont, eux, de dimension finie (voir [Dem93]).

Enfin, les mêmes complexes $\lc^{\bullet}$ et $\bc^{\bullet}$ donnent également une interprétation hypercohomologique de la cohomologie d'Aeppli. On a en effet
$$H^{p,q}_A(X,\C)\cong\Hb^{p+q}(X,\lc_{p+1,q+1}^{\bullet})\cong\Hb^{p+q+1}(X,\bc_{p+1,q+1}^{\bullet}).$$
La cohomologie d'Aeppli est donc elle aussi de dimension finie et calculable à l'aide de formes $C^{\infty}$ ou de courants.

En combinant ces trois remarques, on peut obtenir une preuve plus directe de la dualité entre cohomologie de Bott-Chern de cohomologie d'Aeppli sur une variété compacte : on part du complexe $\lc_{p,q}^{\bullet}$, qui est un complexe d'espaces de Fréchet tel que l'image de la différentielle est fermée (puisque sa cohomologie est de dimension finie, $d\lc^k\subset\ker d$ est de codimension finie et donc fermé). Dès lors le dual de sa cohomologie est donné par la cohomologie du complexe dual. L'espace dual de $\ec^{p,q}(X)$ est l'espace ${\dc'}^{n-p,n-q}$ des courants de bidegré $(n-p,n-q)$, et par exemple le dual de la partie
$$\ldots\to\ec^{p-1,q-1}(X)\stackrel{\ddbar}{\too}\ec^{p,q}(X)\stackrel{d}{\too}\ec^{p+1,q}(X)\oplus\ec^{p,q+1}(X)\to\ldots$$
est donné par
$$\ldots\leftarrow{\dc'}^{n-p+1,n-q+1}(X)\stackrel{\ddbar}{\longleftarrow}{\dc'}^{n-p,n-q}(X)\stackrel{\de+\dbar}{\longleftarrow}{\dc'}^{n-p-1,n-q}(X)\oplus{\dc'}^{n-p,n-q-1}(X)\leftarrow\ldots$$
d'où l'on déduit que le dual de $\lc^{\bullet}_{p,q}$ est le complexe ${\lc'}^{\bullet}_{n-p+1,n-q+1}$ et finalement
$$(H^{p,q}_{BC}(X,\C))^*=(\Hb^{p+q-1}(X,\lc^{\bullet}_{p,q}))^*\cong\Hb^{2n-p-q}(X,{\lc'}^{\bullet}_{n-p+1,n-q+1})=H^{n-p,n-q}_A(X,\C).$$

\subsection{Cohomologie de Bott-Chern entière et cohomologie de Deligne}

\begin{defn} On définit les groupes de cohomologie de Bott-Chern entière par
$$H^{p,q}_{BC}(X,\Z)=\Hb^{p+q}(X,\bc_{\Z(p)}^{\bullet}),$$
où $\bc_{\Z(p)}$ est le complexe $\bc$ dans lequel on a remplacé $\C$ par $\Z(p)=(2\pi i)^p\Z$ en degré 0.\end{defn}
\begin{rem} Le choix de privilégier $p$ par rapport à $q$ en considérant $\Z(p)$ est peu satisfaisant ; mais comme l'objectif est l'étude de classes de Chern, où l'on se limite à $p=q$, ce n'est pas gênant.\end{rem}
Outre l'application naturelle $\eps_{BC}:H^{p,q}_{BC}(X,\Z)\to H^{p,q}_{BC}(X,\C)$, la cohomologie de Bott-Chern entière s'envoie dans la cohomologie de Deligne
$$H^{p+q}_{\dc}(X,\Z(p))=\Hb^{p+q}(X,\Z(p)\to\oc\to\Omega^1\to\ldots\to\Omega^{p-1}\to 0)$$
en ``oubliant'' la partie antiholomorphe, et dans $H^{p+q}_{\dc}(X,\overline{\Z(q)})$ en ``oubliant'' la partie holomorphe et multipliant par $(2\pi i)^{q-p}$ en degré 0. Dans toute la suite, ces applications seront notées :
$$\eps_{\dc}:H^{p,q}_{BC}(X,\Z)\to H^{p+q}_{\dc}(X,\Z(p)),\quad \eps_{\overline{\dc}}:H^{p,q}_{BC}(X,\Z)\to H^{p+q}_{\dc}(X,\overline{\Z(q)}).$$
On remarque que si $q=0$ alors $\eps_{\dc}=\id$. 

\section{Explicitation en cohomologie de \v Cech}

On souhaite expliciter l'isomorphisme entre $H^{p,q}_{BC}(X,\C)$ et $\Hb^{p+q}(X,\bc)$ en cohomologie de \v Cech, ce qui se fait naturellement en utilisant le lemme de résolubilité locale.

\subsection{Hypercohomologie de \v Cech}
Soit $0\to\fc^0\stackrel{d}{\to}\fc^1\to\ldots$ un complexe de faisceaux sur $X$, borné à gauche, et $\Ub$ un recouvrement de $X$. L'hypercohomologie de \v Cech $\check{\Hb}^{\bullet}(\Ub,\fc^{\bullet})$ est la cohomologie du complexe diagonal associé au double complexe $\gc^{p,q}=\check{C}^q(\Ub,\fc^p)$, les différentielles $d$ et $\dcech$ commutant entre elles. Il est clair que pour un recouvrement convenable (par exemple à intersections convexes), $\check{\Hb}^k(\Ub,\fc^{\bullet})\cong\Hb^k(X,\fc^{\bullet})$.

En pratique : une $k$-hypercochaîne de \v Cech de ce complexe est la donnée, pour $0\les j\les k$, d'une $(k-j)$-cochaîne de \v Cech de $\fc^j$. La différentielle de \v Cech est la suivante :
$$\dcechg:(\alpha^0_{i_0...i_k},\alpha^1_{i_0...i_{k-1}},\ldots,\alpha^k_{i_0})\mapsto(\dcech\alpha^0_{i_0...i_k},d\alpha^0_{i_0...i_k}-\dcech\alpha^1_{i_0...i_{k-1}},\ldots,d\alpha^{k-1}_{i_0i_1}+(-1)^k\dcech\alpha^k_{i_0},d\alpha^k_{i_0})$$

\subsection{Isomorphisme entre les hypercohomologies de $\lc^{\bullet}[1]$, $\sca^{\bullet}[1]$ et $\bc^{\bullet}$}$\;$\\
Dans toute la suite on calculera l'hypercohomologie de \v Cech par rapport à un recouvrement $\Ub$ à intersections convexes fixé, et on identifiera $\check{\Hb}^k(\Ub,\fc^{\bullet})$ et $\Hb^k(X,\fc^{\bullet})$.
Un élément de $\Hb^{p+q}(X,\lc^{\bullet}[1])$, avec $p\ges 1$ et $q\ges 1$, est donné par une famille de cochaînes
$$\gamma^{p,q}\in\check{C}^0(\ec^{p,q}),\quad \gamma^{r,s}\in\check{C}^{p+q-r-s-1}(\ec^{r,s}), \; 0\les r\les p-1,\,0\les s\les q-1$$
avec les relations de cocycle
$$d\gamma^{p,q}=0,\qquad \ddbar\gamma^{p-1,q-1}+(-1)^{p+q}\dcech\gamma^{p,q}=0$$
$$\forall 1\les r\les p-1,\;1\les s\les q-1,\quad \de\gamma^{r-1,s}+\dbar\gamma^{r,s-1}+(-1)^{r+s+1}\dcech\gamma^{r,s}=0$$
$$\forall 1\les s\les q-1,\quad \dbar\gamma^{0,s-1}+(-1)^{s+1}\dcech\gamma^{0,s}=0,\qquad \forall 1\les r\les p-1,\quad \de\gamma^{r-1,0}+(-1)^{r+1}\dcech\gamma^{r,0}=0$$
$$\dcech\gamma^{0,0}=0$$

On construit l'élément de $\Hb^{p+q}(X,\sca^{\bullet}[1])$ correspondant en suivant la feuille de route que constitue le lemme de résolubilité locale. On va construire une $(p+q-1)$-hypercochaîne $\alpha=(\alpha^{r,s}_{0\les r\les p-1,0\les s\les q-1})$ de $\lc^{\bullet}$ telle que $\gamma-\dcech\alpha$ soit un hypercocycle de $\sca^{\bullet}$.

\begin{notation} Pour $k\les p+q-2$ donné, on note $r_{\min}=\max(0,k-q+1)$, $r_{\max}=\min(k,p-1)$ et $I_k$ l'ensemble des couples $(r,s)$ avec $r+s=k$ et $r_{\min}\les r\les r_{\max}$ (de sorte que $\lc^k=\bigoplus_{(r,s)\in I_k}\ec^{r,s}$).
\end{notation} 

On va construire $\alpha^{r,s}$ pour $(r,s)\in I_k$, par récurrence descendante sur $k$ :\\
$\diamond$ Pour $k=p+q-2$, $I_k=\{(p-1,q-1)\}$. Comme $\gamma^{p,q}$ est $d$-fermée on peut écrire $\gamma^{p,q}=\ddbar\alpha^{p-1,q-1}$.\\
$\diamond$ Pour $k=p+q-3$ : $\ddbar(\gamma^{p-1,q-1}+(-1)^{p-1+q-1}\dcech\alpha^{p-1,q-1})=0$ donc, si $p\ges 2$ et $q\ges 2$,
$$\gamma^{p-1,q-1}+(-1)^{p+q}\dcech\alpha^{p-1,q-1}=\de\alpha^{p-2,q-1}+\dbar\alpha^{p-1,q-2},$$
et l'on remarque que
$$\de(\gamma^{p-2,q-1}+(-1)^{p-2+q-1}\dcech\alpha^{p-2,q-1})+\dbar(\gamma^{p-1,q-2}+(-1)^{p-1+q-2}\dcech\alpha^{p-1,q-2})=0.$$
Si $p=1$ (resp. $q=1$), on remplace $\de\alpha^{p-2,q-1}$ par $v^{q-1}$ antiholomorphe (resp. $\dbar\alpha^{p-1,q-2}$ par $u^{p-1}$ holomorphe).\\
$\diamond$ Supposons que pour $1\les k\les p+q-3$, on ait défini, pour $(r,s)\in I_k$, les cochaînes $\alpha^{r,s}$, si $k\les r+1$ une cochaîne holomorphe $u^{r+1,0}$, et si $k\les s+1$ une cochaîne antiholomorphe $v^{0,s+1}$, tels que
$$\forall(r,s)\in I_k\moins\{(r_{\max},k-r_{\max})\}, \de(\gamma^{r,s}+(-1)^{r+s}\dcech\alpha^{r,s})+\dbar(\gamma^{r+1,s-1}+(-1)^{r+s}\dcech\alpha^{r+1,s-1})=0.$$ 
On déduit de cette famille de relations une famille de cochaînes $\alpha^{r,s}$ pour $(r,s)\in I_{k-1}$, $r\not=0$ et $s\not=0$, telles que
$$\gamma^{r,s}+(-1)^{r+s}\alpha^{r,s}=\dbar\alpha^{r,s-1}+\de\alpha^{r-1,s}.$$
Si $(0,k)$ (resp. $(k,0)$) apparteint à $I_{k}$, on a à la place
$$\gamma^{0,k}+(-1)^k\dcech\alpha^{0,k}=\dbar\alpha^{0,k-1}+v^{0,k},\mbox{ (resp. } \gamma^{k,0}+(-1)^k\alpha^{k,0}=u^{k,0}+\de\alpha^{k-1,0})$$
avec $v^{0,k}$ antiholomorphe (resp. $u^{k,0}$ holomorphe). On a alors bien, quand elles ont un sens, les relations
$$\de(\gamma^{r,s}+(-1)^{r+s}\dcech\alpha^{r,s})+\dbar(\gamma^{r+1,s-1}+(-1)^{r+s}\dcech\alpha^{r+1,s-1})=0,$$
ce qui permet de poursuivre la récurrence jusqu'à $k=1$ inclus.\\
$\diamond$ Lors de cette dernière étape on a défini $\alpha^{0,0}$, $v^1$ et $u_1$ tels que
$$\gamma^{0,1}-\dcech\alpha^{0,1}=\dbar\alpha^{0,0}+v^{0,1}\mbox{ et }\gamma^{1,0}-\dcech\alpha^{1,0}=u^{1,0}+\de\alpha^{0,0}.$$
$\diamond$ Finalement $\ddbar(\gamma^{0,0}+\dcech\alpha^{0,0})=0$ donc $w^{0,0}=\gamma^{0,0}+\beta^{0,0}$ est une cochaîne de $\oc+\ocbar$.

On a ainsi défini une $(p+q)$-hypercochaîne de $\sca^{\bullet}[1]$ par 
$$w=(w^{0,0},(u^{r,0})_{1\les r\les p-1},(v^{0,s})_{1\les s\les q-1}),$$
et par construction on a bien $\gamma-w=\dcechg\alpha$. Enfin, $w$ n'est modifié que par un cobord si l'on modifie $\gamma$ par un cobord ou si l'on change le choix des $\alpha$.

L'élément correspondant dans $\Hb^{p+q}(X,\bc^{\bullet})$ s'obtient comme suit : on ne change pas $u^{r,0}$ et $v^{0,s}$, on écrit $w^{0,0}=u^{0,0}+v^{0,0}$ avec $u^{0,0}$ holomorphe et $v^{0,0}$ antiholomorphe, et on pose $c=\dcech(u^{0,0})=-\dcech(v^{0,0})$. L'hypercocycle est alors la famille $(c,(u^{r,0})_{0\les r\les p-1},(v^{0,s})_{0\les s\les q-1})$. Une écriture différente $w^{0,0}=\tilde{u}^{0,0}+\tilde{v}^{0,0}$ modifierait cet hypercocycle par l'hypercobord $\dcechg(u^{0,0}-\tilde{u}^{0,0})$.

\subsection{Cas où $p$ ou $q$ est nul}
Supposons $q=0$, l'autre cas étant symétrique. Cette fois un élément de $\Hb^p(X,\lc^{\bullet}[1])$ est simplement donné par un $0$-cocycle $\theta^{p,0}$ tel que $d\theta^{p,0}=0$ ($\theta^{p,0}$ est donc holomorphe $\de$-fermé). L'élément correspondant de $\Hb^p(X,\bc^{\bullet})$ est donné par l'hypercocycle $(c,u^{0,0},\ldots,u^{p-1,0})$ défini par récurrence descendante comme suit :\\
Tout d'abord $\theta^{p,0}$ est holomorphe $\de$-fermé donc $\theta^{p,0}=\de u^{p-1,0}$ avec $u^{p-1,0}$ holomorphe, et $\de\dcech u^{p-1,0}=0$.\\
Supposons $u^{r,0}$ construit $(r\ges 1)$ tel que $\de\dcech u^{r,0}=0$, alors $(-1)^r\dcech u^{r,0}=\de u^{r-1,0}$, avec $u^{r-1,0}$ holomorphe, et $\de\dcech u^{r-1,0}=0$. Finalement $\dcech u^{0,0}=c$.

\subsection{Bilan}
Soit un élément de $H^{p,q}_{BC}(X,\C)$, représenté par une $(p,q)$-forme fermée $\theta$. Il est défini dans $\Hb^{p+q}(X,\lc[1]^{\bullet})$ par l'hypercocycle, encore notée $\theta$ et définie par $\theta^{p,q}=\theta|_{U_j}$ et $\theta^{r,s}=0$ sinon. 
On applique à $\theta$ la construction précédente : si $p\ges 1$ et $q\ges 1$, il existe un hypercocycle $w=(c;u^{r,0};v^{0,s})\in\check{Z}^{p+q}(X,\bc^{\bullet})$ et une hypercochaîne $\alpha=(\alpha^{r,s})\in\check{C}^{p+q-1}(X,\lc[1]^{\bullet})$ tels que $\theta=\dcechg\alpha+w$. On représentera ces données sous forme du tableau suivant :
$$\theta\longleftrightarrow\left[
\begin{array}{c|ccc}
v^{0,q-1}&&&\\
\vdots&&\alpha^{r,s}&\\
v^{0,0}&&&\\\hline
c&u^{0,0}&\cdots&u^{p-1,0}\end{array}\right]$$
L'égalité $\theta=\dcechg\alpha+w$ correspond aux relations suivantes :
$$(\bigstar)\left\{\begin{array}{rcll}
\theta^{p,q}&=&\ddbar\alpha^{p-1,q-1}&\\
(-1)^{r+s}\dcech\alpha^{r,s}&=&\dbar\alpha^{r,s-1}+\de\alpha^{r-1,s}&\forall\, 1\les r\les p-1,\, 1\les s\les q-1\\
(-1)^s\dcech\alpha^{0,s}&=&\dbar\alpha^{0,s-1}+v^{0,s}&\forall\, 1\les s\les q-1\\
(-1)^r\dcech\alpha^{r,0}&=&u^{r,0}+\de\alpha^{r-1,0}&\forall\, 1\les r\les p-1\\
\dcech\alpha^{0,0}&=&u^{0,0}+v^{0,0}&\\
\dcech u^{0,0}&=&c&
\end{array}\right.$$
Remarquons que ces relations impliquent les relations d'hypercocycle pour $u$ et $v$ :
$$(-1)^r\dcech u^{r,0}=\de u^{r-1,0}\;\forall 1\les r\les p-1,\qquad (-1)^s\dcech u^{0,s}=\dbar v^{0,s-1}\;\forall 1\les s\les q-1$$

Si $q=0$, on a simplement
$$\theta\longleftrightarrow\left(c,u^{0,0},\ldots,u^{p-1,0}\right)$$
avec les relations
$$\theta^{p,0}=\de u^{p-1,0},\quad (-1)^r\dcech u^{r,0}=\de u^{r-1,0}\;\forall 1\les r\les p-1,\quad \dcech u^{0,0}=c.$$
De même si $p=0$, on a 
$$\theta\longleftrightarrow\left(c,v^{0,0},\ldots,v^{0,q-1}\right)$$
avec les relations
$$\theta^{0,q}=-\dbar v^{0,q-1},\quad (-1)^s\dcech v^{0,s}=\dbar v^{0,s-1}\;\forall 1\les r\les p-1,\quad -\dcech v^{0,0}=c.$$

\section{Structure d'algèbre sur les groupes de cohomologie de Bott-Chern}

Le but de cette partie est de munir la cohomologie de Bott-Chern entière d'une structure d'algèbre, c'est-à-dire de définir un produit
$$H^{p,q}_{BC}(X,\Z)\times H^{p',q'}_{BC}(X,\Z)\too H^{p+p',q+q'}_{BC}(X,\Z).$$
Il se trouve que pour la cohomologie de Bott-Chern classique, nous disposons d'une telle structure, donnée par le produit extérieur. Il s'agit donc de voir comment le produit extérieur ``traverse" l'isomorphisme $H^{p,q}_{BC}(X,\C)\cong\Hb^{p+q}(X,\bc)$, et d'adopter la définition ainsi obtenue.

\subsection{Cup-produit}

On continue à travailler en cohomologie de \v Cech, à l'aide du cup-produit. Plus précisément, si $\beta^{p,q}\in\check{C}^k(\ec^{p,q})$ et $\tilde{\beta}^{p',q'}\in\check{C}^{\ell}(\ec^{p',q'})$, on combine cup-produit et produit extérieur pour définir $\beta^{p,q}\cdot\tilde{\beta}^{p',q'}\in\check{C}^{k+\ell}(\ec^{p+p',q+q'})$ par
$$(\beta\cdot\tilde{\beta})_{j_0\ldots j_{k+\ell}}=\beta_{j_0\ldots j_k}\wedge\tilde{\beta}_{j_k\ldots j_{k+\ell}}.$$
On a les relations 
$$\dcech(\beta^{p,q}\cdot\tilde{\beta}^{p',q'})=(\dcech\beta^{p,q})\cdot\tilde{\beta}^{p',q'}+(-1)^k\beta^{p,q}\cdot(\dcech\tilde{\beta}^{p',q'})$$
$$\de(\beta^{p,q}\cdot\tilde{\beta}^{p',q'})=(\de\beta^{p,q})\cdot\tilde{\beta}^{p',q'}+(-1)^{p+q}\beta^{p,q}\cdot(\de\tilde{\beta}^{p',q'})$$
$$\dbar(\beta^{p,q}\cdot\tilde{\beta}^{p',q'})=(\dbar\beta^{p,q})\cdot\tilde{\beta}^{p',q'}+(-1)^{p+q}\beta^{p,q}\cdot(\dbar\tilde{\beta}^{p',q'})$$

\subsection{Traduction du produit extérieur}
\begin{prop} On se donne $p,p',q,q'$ tels qu'on n'ait pas $p'=0$ et $p\not=0$, ou $q=0$ et $q'\not=0$. Soient $\{\theta\}\in H^{p,q}_{BC}(X,\C)$ et $\{\tilde{\theta}\}\in H^{p',q'}_{BC}(X,\C)$, donnés par les hypercocycles 
$$w=\left(c,\begin{array}{l} u^{0,0},\ldots, u^{p-1,0}\\ v^{0,0},\ldots\ldots, v^{0,q-1}\end{array}\right),\qquad \tilde{w}=\left(\tilde{c},\begin{array}{l} \tilde{u}^{0,0},\ldots\ldots\ldots, \tilde{u}^{p'-1,0}\\ \tilde{v}^{0,0},\ldots, \tilde{v}^{0,q'-1}\end{array}\right).$$
Alors $\{\theta\wedge\tilde{\theta}\}$ est donné \textbf{aux signes près} (que l'on précisera) par l'hypercocycle
$$W=\left(c\cdot\tilde{c},\begin{array}{l} c\cdot\tilde{u}^{0,0}\,,\,\ldots\ldots\ldots,c\cdot\tilde{u}^{p'-1,0}\,,\,u^{0,0}\cdot\de\tilde{u}^{p'-1,0}\,,\,\ldots\,,\,u^{p-1,0}\cdot\de\tilde{u}^{p'-1,0}\\ v^{0,0}\cdot\tilde{c}\,,\,\ldots\,,\,v^{0,q-1}\cdot\tilde{c}\,,\,\dbar v^{0,q-1}\cdot\tilde{v}^{0,0}\,,\,\ldots\ldots\,,\,\dbar v^{0,q-1}\cdot\tilde{v}^{0,q'-1}\end{array}\right).$$
\end{prop}
\begin{proof} 
On écrit $\theta=\dcechg\alpha+w$ et $\tilde{\theta}=\dcechg\tilde{\alpha}+\tilde{w}$. Un calcul direct montre alors que $\theta\wedge\tilde{\theta}=\dcechg A+W$ avec :
$$A^{R,S}=\left\{\begin{array}{clll}
\eps^{R,S}\dbar\alpha^{R-p',q-1}\cdot\de\tilde{\alpha}^{p'-1,S-q}&\mbox{si}&p'\les R\les p+p'-1,& q\les S\les q+q'-1\\
\eps^{R,S}\dbar v^{0,q-1}\cdot\tilde{\alpha}^{R,S-q}&\mbox{si}&0\les R\les p'-1,& q\les S\les q+q'-1\\
\eps^{R,S}\alpha^{R-p',R}\cdot\de\tilde{u}^{p'-1,0}&\mbox{si}&p'\les R\les p+p'-1,& 0\les S\les q-1\\
\eps^{R,S}v^{0,S}\cdot\de\tilde{u}^{R,0}&\mbox{si}&0\les R\les p'-1,& 0\les S\les q-1\\
\end{array}\right.$$
et $W=(C,U^{R,0},V^{0,S})$ comme dans l'énoncé, à savoir
$$U^{R,0}=\left\{\begin{array}{cll}
\eps^{R,\star}u^{R-p',0}\cdot\de\tilde{u}^{p'-1,0}&\mbox{si}&p'\les R\les p+p'-1\\
\eps^{R,\star}c\cdot\de\tilde{u}^{R,0}&\mbox{si}&0\les R\les p'-1\\
\end{array}\right.$$
$$V^{0,S}=\left\{\begin{array}{cll}
\eps^{\star,S}\dbar v^{0,q-1}\cdot\tilde{v}^{0,S-q}&\mbox{si}&q\les S\les q+q'-1\\
\eps^{\star,S}v^{0,S}\cdot\tilde{c}&\mbox{si}& 0\les S\les q-1\\
\end{array}\right.$$
$$C=\eps^{\star,\star}c\cdot\tilde{c},$$
où tous les symboles $\eps^{\bullet,\bullet}$ sont des signes que l'on va déterminer

Cette expression se lit plus agréablement à l'aide de l'écriture en tableaux :
$$\theta\longleftrightarrow\left[
\begin{array}{c|ccc}
v^{0,q-1}&&&\\
\vdots&&\alpha^{r,s}&\\
v^{0,0}&&&\\\hline
c&u^{0,0}&\cdots&u^{p-1,0}\end{array}\right],
\quad\tilde{\theta}\longleftrightarrow\left[
\begin{array}{c|ccc}
\tilde{v}^{0,q'-1}&&&\\
\vdots&&\tilde{\alpha}^{r',s'}&\\
\tilde{v}^{0,0}&&&\\\hline
\tilde{c}&\tilde{u}^{0,0}&\cdots&\tilde{u}^{p-1,0}\end{array}\right]$$

$$\theta\wedge\tilde{\theta}\longleftrightarrow\left[
\begin{array}{c|cccc}
\eps^{\star,S}\dbar v^{0,q-1}\cdot\tilde{v}^{0,S-q}&&\eps^{R,S}\dbar v^{0,q-1}\cdot\tilde{\alpha}^{R,S-q}&&\eps^{R,S}\dbar\alpha^{R-p',q-1}\cdot\de\tilde{\alpha}^{p'-1,S-q}\\
\\
\eps^{\star,S}v^{0,S}\cdot\tilde{c}&&\eps^{R,S}v^{0,S}\cdot \tilde{u}^{R,0}&& \eps^{R,S}\alpha^{R-p',S}\cdot\de\tilde{u}^{p'-1,0}\\
\\\hline
\eps^{\star,\star}c\cdot\tilde{c}&&\eps^{R,\star}c\cdot\tilde{u}^{R,0}&&\eps^{R,\star}u^{R-p',0}\cdot\de\tilde{u}^{p'-1,0}
\end{array}\right]$$

Pour obtenir les signes, on substitue les expressions ci-dessus dans les relations $(\bigstar)$ pour $\theta\wedge\tilde{\theta}$. On obtient les relations de récurrence suivantes :\\
Abaissement du deuxième indice, pour $p'\les R\les p+p'-1$ :
$$\begin{array}{rcl}
\eps^{R,S-1}&=&\left\{ \begin{array}{ll} (-1)^{p+p'+R+1}\eps^{R,S}&\mbox{si }S\ges q+1 \\ (-1)^{p+q}\eps^{R,S}&\mbox{si }S=q \\ (-1)^{p'}\eps^{R,S}&\mbox{si }S\les q-1\end{array}  \right. \\
\eps^{R,\star}&=&(-1)^{p'}\eps^{R,0}
\end{array}$$
Abaissement du deuxième indice, pour $0\les R\les p'-1$ :
$$\begin{array}{rcl}
\eps^{R,S-1}&=&\left\{ \begin{array}{ll} (-1)^p\eps^{R,S}&\mbox{si }S\ges q+1 \\ (-1)^{p+q}\eps^{R,S}&\mbox{si }S=q \\ (-1)^R\eps^{R,S}&\mbox{si }S\les q-1\end{array}  \right. \\
\eps^{R,\star}&=&(-1)^{R+1}\eps^{R,0}
\end{array}$$
Abaissement du deuxième indice, pour $R=\star$ :
$$\begin{array}{rcl}
\eps^{\star,S-1}&=&\left\{ \begin{array}{ll} (-1)^p\eps^{\star,S}&\mbox{si }S\ges q+1 \\ (-1)^{p+q+1}\eps^{\star,S}&\mbox{si }S=q \\ \eps^{\star,S}&\mbox{si }S\les q-1\end{array}  \right. \\
\eps^{\star,\star}&=&\eps^{\star,0}
\end{array}$$
Abaissement du premier indice, pour $q\les S\les q+q'-1$ :
$$\begin{array}{rcl}
\eps^{R-1,S}&=&\left\{ \begin{array}{ll} (-1)^{p'+q+S}\eps^{R,S}&\mbox{si }R\ges p'+1 \\ (-1)^{p'+S+1}\eps^{R,S}&\mbox{si }R=p' \\ (-1)^p\eps^{R,S}&\mbox{si }R\les p'-1\end{array}  \right. \\
\eps^{\star,S}&=&(-1)^{p+q}\eps^{0,S}
\end{array}$$
Abaissement du premier indice, pour $0\les S\les q-1$ :
$$\begin{array}{rcl}
\eps^{R-1,S}&=&\left\{ \begin{array}{ll} (-1)^{p'}\eps^{R,S}&\mbox{si }R\ges p'+1 \\ (-1)^{p'+S}\eps^{R,S}&\mbox{si }R=p' \\ (-1)^{p+q+S+1}\eps^{R,S}&\mbox{si }R\les p'-1\end{array}  \right. \\
\eps^{\star,S}&=&(-1)^{p+q+1}\eps^{0,S}
\end{array}$$
Abaissement du premier indice, pour $S=\star$ :
$$\begin{array}{rcl}
\eps^{R-1,\star}&=&\left\{ \begin{array}{ll} (-1)^{p'}\eps^{R,\star}&\mbox{si }R\ges p'+1 \\ (-1)^{p'}\eps^{R,\star}&\mbox{si }R=p' \\ (-1)^{p+q}\eps^{R,\star}&\mbox{si }R\les p'-1\end{array}  \right. \\
\eps^{\star,\star}&=&(-1)^{p+q}\eps^{0,\star}
\end{array}$$

Pour initialiser la récurrence, il faut considérer différents cas :\\
Si $p\ges 1$, $p'\ges 1$, $q\ges 1$ et $q'\ges 1$ alors $\eps^{p+p'-1,q+q'-1}=(-1)^{p+q}$ d'après
$$\ddbar(\alpha^{p-1,q-1})\cdot\ddbar(\tilde{\alpha}^{p'-1,q'-1})=(-1)^{p+q}\ddbar(\dbar\alpha^{p-1,q-1}\cdot\de\tilde{\alpha}^{p'-1,q'-1}).$$
Si $p\ges 1$, $p'\ges 1$, $q\ges 1$ et $q'=0$ alors $\eps^{p+p'-1,q-1}=1$ d'après
$$\ddbar(\alpha^{p-1,q-1})\cdot\de(\tilde{u}^{p'-1,0})=\ddbar(\alpha^{p-1,q-1}\cdot\de\tilde{u}^{p'-1,0}).$$
Si $p=0$, $p'\ges 1$, $q\ges 1$ et $q'\ges 1$ alors $\eps^{p'-1,q+q'-1}=1$ d'après
$$-\dbar v^{0,q-1}\cdot\ddbar(\tilde{\alpha}^{p'-1,q'-1})=\ddbar(\dbar v^{0,q-1}\cdot\tilde{\alpha}^{p'-1,q'-1}).$$
Si $p=0$, $p'\ges 1$, $q\ges 1$ et $q'=0$ alors $\eps^{p'-1,q+-1}=(-1)^{q}$ d'après
$$-\dbar(v^{0,q-1})\cdot\de(\tilde{u}^{p'-1,0})=(-1)^q\ddbar(v^{0,q-1}\cdot\de\tilde{u}^{p'-1,0}).$$
Si $p\ges 1$, $p'\ges 1$, $q=q'=0$ alors $\eps^{p+p'-1,\star}=1$ d'après
$$\de u^{p-1,0}\cdot\de\tilde{u}^{p'-1,0}=\de(u^{p-1,0}\cdot\de\tilde{u}^{p'-1,0}).$$
Si $p=p'=0$, $q\ges 1$ et $q'\ges 1$ alors $\eps^{\star,q+q'-1}=(-1)^{q+1}$ d'après
$$(-\dbar v^{0,q-1})\cdot(-\dbar\tilde{v}^{0,q'-1})=-(-1)^{q-1}\dbar(\dbar v^{0,q-1}\cdot\tilde{v}^{0,q'-1}).$$

On constate que ces différentes initialisations donnent en fin de compte les mêmes formules. Voici le tableau donnant le signe $\eps^{R,S}$ :

\renewcommand{\arraystretch}{1.2}
$$\begin{array}{r|c|c|c|}\hline
S\ges q & (-1)^{pS+(p+1)(q+1)} & (-1)^{p(R+S+q)+1} & \displaystyle{(-1)^{(R+p+p'+1)(S+q+p')+p+q}} \\\hline
S\les q-1 & 1 & (-1)^{R(S+p+q+1)+p+q+1} & (-1)^{p'(R+S+p+q+1)} \\\hline
S=\star & 1 & (-1)^{(p+q)(R+1)} & (-1)^{p'(R+p+q)} \\\hline
& R=\star & R\les p'-1 & R\ges p'
\end{array}$$
\renewcommand{\arraystretch}{1}

\end{proof}

Comme annoncé, on utilise les formules ci-dessus pour mettre une structure d'anneau sur $H^{\bullet,\bullet}_{BC}(X,\Z)$.

\begin{rems} 1. Le produit extérieur étant anti-commutatif, le produit ainsi défini l'est aussi (cela pourrait se vérifier directement). Ceci permet de traiter les cas $p\not=0$, $p'=0$ et $q=0$, $q'\not=0$.

2. La formule obtenue est compatible, aux signes près, avec le produit en cohomologie de Deligne, tel qu'il est par exemple défini dans [EV88].
\end{rems}

%

\section{Eléments de structure des groupes de cohomologie de Bott-Chern entière}

On va donner ici quelques indications sur la structure des groupes $H^{p,q}_{BC}(X,\Z)$ à l'aide de suites exactes faisant intervenir des groupes connus.

\subsection{Lien avec la cohomologie de Bott-Chern usuelle}
On a le diagramme commutatif suivant :
$$\begin{array}{ccc} H^{p,q}_{BC}(X,\Z)&\stackrel{\eps_{BC}}{\too}&H^{p,q}_{BC}(X,\C)\\
\downarrow\sm{\fy_\Z}&&\downarrow\sm{\fy}\\
H^{p+q}(X,\Z(p))&\stackrel{\eps}{\too}&H^{p+q}(X,\C)\end{array}$$
\begin{prop}\label{produitfibre} Soient $\xi\in H^{p,q}_{BC}(X,\C)$ et $\alpha\in H^{p+q}(X,\Z(p))$ des classes telles que $\fy(\xi)=\eps(\alpha)$. Alors il existe une classe $\tilde{\xi}\in H^{p,q}_{BC}(X,\Z)$ telle que $\xi=\eps_{BC}(\tilde{\xi})$ et $\alpha=\fy_{\Z}(\tilde{\xi})$. 
\end{prop}
\begin{proof} Soit $\theta=(c;u;v)$ un hypercocycle représentant $\xi$, et $a$ un cocycle représentant $\alpha$. Par hypothèse $\{c\}=\{a\}\in H^{p+q}(X,\C)$, donc il existe une cochaîne $b\in\check{C}^{p+q-1}(X,\C)$ telle que $c-\dcech b=a$. On pose $\tilde{\theta}=\theta-\dcechg(b;0;0)$ qui est un hypercocycle de $\bc_{\Z(p)}$, et dont la classe $\tilde{\xi}\in H^{p,q}_{BC}(X,\Z)$ convient.\end{proof}

\begin{rem} Il n'y a pas en général unicité de cet élément, comme le montrera un peu plus loin l'exemple de l'espace projectif.\end{rem}

Il est également intéressant de considérer la suite exacte courte de complexes
$$0\too\bc_{\Z(p)}\too\bc\too\C/\Z\too0$$
et la suite exacte longue associée qui s'écrit, compte tenu de l'interprétation hypercohomologique de la cohomologie d'Aeppli (voir la remarque du paragraphe 2.c.)
$$H^{p-1,q-1}_A(X,\C)\too H^{p+q-1}(X,\C/\Z)\too H^{p,q}_{BC}(X,\Z)\too H^{p,q}_{BC}(X,\C)$$

\subsection{Analogue de la suite exacte exponentielle}

La suite exacte exponentielle classique est la suite exacte courte donnée par
$$0\too\Z\too\oc\stackrel{\exp(2\pi i\cdot)}{\too}\oc^*\too 0.$$
Une retraduction de l'exactitude de cette suite est le fait que le complexe $\oc^*[1]$, formé de l'unique terme $\oc^*$ en position $1$, est quasi-isomorphe au complexe de Deligne $0\to\Z\to\oc\to 0$ et, sous ce quasi-isomorphisme, la suite exacte exponentielle s'apparente à la suite exacte courte de complexes de faisceaux
$$0\too\oc[1]\too(\Z\to\oc)\too\Z\too 0.$$
On cherche ici un analogue de cette suite exacte.

On définit pour $p\ges 1$ les complexes de De Rham holomorphes tronqués par
$$\begin{array}{rc}
\Omega_{\ges p}^{\bullet}=F^p\Omega^{\bullet} :& 0\to\Omega^{p}\to\Omega^{p+1}\to\ldots\\
$$\Omega_{<p}^{\bullet} :& \oc\to\Omega^1\to\ldots\to\Omega^{p-1}\to 0\end{array}$$

Supposons que $p\ges 1$ et $q\ges 1$. On considère la suite exacte courte évidente
$$0\too(\Omega^{\bullet}_{<p}\oplus\ombar{\bullet}_{<q})[1]\too\bc_{\Z(p)}\too\Z(p)\too 0$$
et la suite exacte longue d'hypercohomologie correspondante, qui fait intervenir les groupes $\Hb^k(X,\Omega_{<p}^{\bullet})$. 
Compte tenu de l'isomorphisme entre $\Z$ et $\Z(p)$, on obtient :
$$H^{p+q-1}(X,\Z)\to\Hb^{p+q-1}(X,\Omega_{<p}^{\bullet}\oplus\overline{\Omega}^{\bullet}_{<q})\to H^{p,q}_{BC}(X,\Z)\to H^{p+q}(X,\Z)\to\Hb^{p+q}(X,\Omega_{<p}^{\bullet}\oplus\overline{\Omega}^{\bullet}_{<q}).$$

Ces groupes peuvent être calculés, dans le cas kählérien compact, à l'aide de la décomposition de Hodge :

\begin{lemma}\label{kahler} Supposons $X$ kählérienne compacte. Alors pour $p\in\N^*$ et $k\in\N$ on a
$$\Hb^k(X,\Omega_{<p}^{\bullet})=\bigoplus_{\substack{r+s=k\\ r< p}}H^{r,s}(X,\C).$$
\end{lemma}
\begin{proof}
Il ressort comme conséquence de la théorie de Hodge que le morphisme
$$\Hb^k(X,\Omega_{\ges p}^{\bullet})\to\Hb^k(X,\Omega^{\bullet})=H^k(X,\C)$$
est injectif, d'image
$$\bigoplus_{\substack{r+s=k\\ r\ges p}}H^{r,s}(X,\C).$$
On conclut à l'aide de la suite exacte courte  
$0\to\Omega_{\ges p}^{\bullet}\to\Omega^{\bullet}\to\Omega_{<p}^{\bullet}\to 0$.
\end{proof}

En particulier
$$\Hb^{p+q-1}(X,\Omega_{<p}^{\bullet}\oplus\overline{\Omega}^{\bullet}_{<q})=\bigoplus_{r+s=p+q-1}H^{r,s}(X,\C)=H^{p+q-1}(X,\C)$$
$$\Hb^{p+q}(X,\Omega_{<p}^{\bullet}\oplus\overline{\Omega}^{\bullet}_{<q})=\bigoplus_{\substack{r+s=p+q\\(r,s)\not=(p,q)}}H^{r,s}(X,\C)=H^{p+q}(X,\C)/H^{p,q}(X,\C),$$
et la suite exacte longue précédente se réécrit
$$H^{p+q-1}(X,\Z)\too H^{p+q-1}(X,\C)\too H^{p,q}_{BC}(X,\Z)\too H^{p+q}(X,\Z)\too H^{p+q}(X,\C)/H^{p,q}(X,\C).$$

\subsection{Groupes $H^{p,p}_{BC}(X,\Z)$ et cohomologie de Deligne}

On considère cette fois la suite exacte courte
$$0\too(\ombar{\bullet}_{<q})[1]\too\bc_{\Z(p)}\too\Z(p)_{\dc}\too 0$$ 
qui relie la cohomologie de Bott-Chern entière à la cohomologie de Deligne. La situation est particulièrement intéressante dans le cas où $p=q$ :

\begin{prop}\label{BCZDeligne} Soit $p\ges 1$. Alors
$$H^{p,p}_{BC}(X,\Z)\cong H^{2p}_{\dc}(X,\Z(p))\oplus \overline{\Hb^{2p-1}(X,\Omega_{<p}^{\bullet})}.$$
En particulier, si $X$ est kählerienne compacte,
$$H^{p,p}_{BC}(X,\Z)\cong H^{2p}_{\dc}(X,\Z(p))\oplus\bigoplus_{\substack{r+s=2p-1\\ r>s}}H^{r,s}(X,\C).$$
\end{prop}
\begin{proof} On remarque simplement que, dans le cas $p=q$, la suite exacte courte précédente est scindée par le morphisme 
$$\begin{array}{rcl}\Omega^{k}&\too&\Omega^k\oplus\ombar{k}\\
u&\longmapsto &\left(u,(-1)^{p+1}\overline{u}\right),\end{array}$$
Le signe $(-1)^{p+1}$ est nécessaire pour rendre le diagramme
$$\begin{array}{ccl}
(2\pi i)^p\Z&\stackrel{(+,-)}{\too}&\oc\oplus\ocbar\\
\downarrow\uparrow =&&\quad\downarrow\uparrow\sm{u,(-1)^{p+1}\overline{u}}\\
(2\pi i)^p\Z&\too&\quad\oc\end{array}$$
commutatif, et la commutativité se propage jusqu'à la fin du complexe, y compris au dernier rang car $p=q$.

Le cas kählerien provient du lemme \ref{kahler}.
\end{proof}

\subsection{Cohomologie de l'espace projectif}

On va appliquer ces résultats pour calculer la cohomologie de l'espace projectif.
\begin {prop} 
1. Pour $p\in\{0,\ldots,n\}$ on a $H^{p,p}_{BC}(\Pb^n,\Z)\cong \Z$, et on a l'isomorphisme d'algèbre
$$\bigoplus_pH^{p,p}_{BC}(\Pb^n,\Z)\cong\Z[h]/h^{n+1}$$
où $h$ est un générateur de $H^{1,1}_{BC}(\Pb^n,\Z)$.\\
2. Si $p\not=q$ alors
$$H^{p,q}_{BC}(X,\Z)\cong\left\{ \begin{array}{ll} 0&\mbox{si }p+q\mbox{ est pair}\\\C/\Z&\mbox{si }p+q\mbox{ est impair}\end{array}  \right.$$
\end{prop}
\begin{proof} Le premier point provient de la proposition \ref{BCZDeligne} : le terme 
$$\bigoplus_{\substack{r+s=2p-1\\ r>s}}H^{r,s}(\Pb^n,\C)$$ 
est nul pour l'espace projectif, de sorte que la sous-algèbre $\bigoplus_pH^{p,p}_{BC}(\Pb^n,\Z)$ s'identifie à la même sous-algèbre pour la cohomologie de Deligne, pour laquelle le résultat est vrai.

Dans le cas où $p\not=q$, on considère la suite exacte du paragraphe 5.a., et le fait que si $p\not=q$, alors $H^{p,q}_{BC}(\Pb^n,\C)=0$ et $H^{p-1,q-1}_{A}(\Pb^n,\C)=0$. On a donc
$$H^{p,q}_{BC}(\Pb^n,\Z)\cong H^{p+q-1}(\Pb^n,\C/\Z),$$
d'où le résultat annoncé.\end{proof}

\section{Classes de Chern en cohomologie de Bott-Chern entière}

On souhaite définir des classes de Chern dans ces groupes de cohomologie. En fait, il suffit de le faire dans le cas des fibrés en droites, des techniques classiques permettant de passer au fibrés vectoriels de tout rang puis aux faisceaux cohérent et en particulier aux cycles analytiques.

\subsection{Expression de la classe de Chern d'un fibré en droites dans les différentes cohomologies}
Soit $L$ un fibré en droites holomorphe sur $X$, et $\Ub=(U_j)$ un recouvrement à intersections convexes tel que sur $U_j$, $L$ soit trivialisé par une section $e_j$ partout non nulle. On note $g_{jk}$ la fonction de transition définie sur $U_j\cap U_k$ par $e_k(x)=g_{jk}(x)e_j(x)$, l'élément 
$$\{g_{jk}\}\in\check{H}^1(\Ub,\oc^*)\cong H^1(X,\oc^{\star})$$ 
déterminant la classe d'isomorphisme de $L$.

Le complexe de Deligne $\Z(1)_{\dc}:\Z(1)\to\oc\to 0$ est quasi-isomorphe au complexe $\oc^*[1]$, via l'exponentielle $\exp:\oc\to\oc^*$ dont le noyau est précisément le faisceau $\Z(1)$. Les espaces $H^2_{\dc}(X,\Z(1))$ et $H^1(X,\oc^*)$ sont donc isomorphes, et on note $c_1(L)_{\dc}$ l'image par cet isomorphisme de l'élément $\{g_{jk}\}\in H^1(X,\oc^*)$. Explicitement, quitte à raffiner le recouvrement $\Ub$ on peut supposer que $g_{jk}=\exp(u_{jk})$, et la condition de cocycle pour $g_{jk}$ implique que $\dcech(u_{jk})=(2\pi i c_{jkl})\in\check{Z}^2(X,\Z(1))$. Finalement
$$c_1(L)_{\dc}=\left\{(2\pi i c_{jkl}),(u_{jk})\right\}\in H^2_{\dc}(X,\Z(1)).$$

L'image par le morphisme d'hypercohomologie induit par la projection $\Z(1)_{\dc}\to\Z(1)$ donne, après division par $2\pi i$, la première classe de Chern usuelle
$$c_1(L)_{\Z}=\left\{(c_{jkl})\right\}\in H^2(X,\Z),$$
et, via l'inclusion de faisceaux $\Z\subset\C$, la même formule définit $c_1(L)_{\C}\in H^2(X,\C)$.\\

Munissons à présent $L$ d'une métrique hermitienne $h$. Soit $D$ la connexion de Chern associée à $(L,h)$ et $\Theta$ la courbure de cette connexion. C'est une $(1,1)$-forme $d$-fermée qui définit donc une classe de $H^{1,1}_{BC}(X,\C)$, cette classe étant indépendante de la métrique hermitienne choisie. On définit usuellement la première classe de Chern de $L$ en cohomologie de Bott-Chern comme la classe 
$$c_1(L)_{BC}=\left\{\frac{i}{2\pi}\Theta\right\}\in H^{1,1}_{BC}(X,\C),$$
de sorte que l'image de cette classe dans la cohomologie de De Rham $H^2(X,\C)$ coïncide avec $c_1(L)_{\C}$ via l'isomorphisme de De Rham-Weil. Toutefois, comme $\Z(1)$ apparaît préféren\-tiellement à $\Z$, c'est la classe $2\pi ic_1(L)_{BC}=\{-\Theta\}$ que nous allons considérer.

Sur $U_j$, la connexion de Chern est donnée par
$$D(\xi_j(x)e_j(x))=(d\xi_j(x)-\de\fy_j(x)\xi_j(x))\otimes e_j(x))$$
où $\fy_j$ est le poids de la métrique dans la trivialisation, défini par la formule 
$$e^{-\fy_j(z)}=|e_j(z)|_h^2,$$ 
et vérifiant la condition de compatibilité sur $U_j\cap U_k$ :
$$-\fy_k+\fy_j=u_{jk}+\overline{u_{jk}}.$$
La courbure sur $U_j$ est de ce fait donnée par $\Theta|_{U_j}=\ddbar \fy_j$, et l'hypercocycle de $\bc$ correspondant à $-\Theta$ est donc
$$\{-\Theta\}\longleftrightarrow \left\{(2\pi i c_{jkl}),(u_{jk}),(\overline{u_{jk}})\right\}.$$
Les $c_{jkl}$ étant entiers, on a défini la première classe de Chern de $L$ en cohomologie de Bott-Chern entière :
$$c_1(L)_{BC,\Z}=\left\{(2\pi i c_{jkl}),(u_{jk}),(\overline{u_{jk}})\right\}\in H^{1,1}_{BC}(X,\Z),$$
de sorte que 
$$\eps_{\dc}c_1(L)_{BC,\Z}=c_1(L)_{\dc},\quad \eps_{BC}c_1(L)_{BC,\Z}=2\pi ic_1(L)_{BC}.$$

\subsection{Classes de Chern des fibrés vectoriels et des faisceaux cohérents}

Pour une variété complexe lisse $X$, on souhaite définir pour tout fibré vectoriel holomorphe $E$ de rang $r$ des classes de Chern $c_k(E)_{BC,\Z}\in H^{k,k}_{BC}(X,\Z)$. On utilise la méthode d'espace classifiant et de principe de scindage de A. Grothendieck [Gro56] :
\begin{prop} Les classes de Chern sont uniquement déterminées par les conditions suivantes :\\
Fonctorialité : pour toute application holomorphe $f:Y\to X$ et tout fibré $E$ sur $X$ on a 
$$f^*c_k(E)_{BC,\Z}=c_k(f^*E)_{BC,\Z}\in H^{k,k}_{BC}(Y,\Z);$$
Compatibilité avec les classes de Chern usuelles : pour tout fibré $E$ sur $X$, l'image de $c_k(E)_{BC,\Z}$ dans $H^{2k}(X,\Z)$ coïncide avec la classe de Chern usuelle $c_k(E)$.
\end{prop}
\begin{thm}[Principe de scindage] Soit $E$ un fibré vectoriel de rang $r$ sur la variété $X$, alors il existe une variété $Y$ et une application holomorphe $f:Y\to X$ telles que :\\
1. Le morphisme de cohomologie $f^*:H^{\bullet,\bullet}_{BC}(X,\Z)\to H^{\bullet,\bullet}_{BC}(Y,\Z)$ est injectif,\\
2. Le fibré image réciproque $f^*E$ sur $Y$ est filtré par des fibrés en $0=E_0\subset E_1\subset\cdots\subset E_r=f^*E$ sur $Y$ dont les quotients successifs sont des fibrés en droites.\end{thm}
\begin{proof} On prend pour $Y$ la variété des drapeaux complets
$$v=\{0=V_0\subset V_1\subset\ldots\subset V_r=V\}$$
pour chaque fibre $V=E_x$, avec la projection naturelle $f:Y\to X$, et on note $E_j$ le fibré universel tel que $E_{j,v}=V_j$. Le point 2. est alors bien vérifié. En fait cette variété de drapeaux se réalise de proche en proche en considérant tout d'abord le fibré en espaces projectifs $\Pb(E)\to X$ [...]
Localement, $\Pb(E)|_U$ a la structure du produit $U\times \Pb^{r-1}$. On souhaite appliquer la formule de Künneth pour justifier que 
\end{proof}
\begin{defn} Avec les notations précédentes, la classe de Chern de $E$ est l'unique élément
$$c_{\bullet}(E)\in\bigoplus_{k=0}^nH^{k,k}_{BC}(X,\Z)$$
tel que
$$f^*c_{\bullet}(E)=\prod_{j=1}^r(1+c_1(E_{j}/E_{j-1})_{BC,\Z})\in\bigoplus H^{k,k}_{BC}(Y,\Z)$$\end{defn}

On souhaite maintenant définir les classes de Chern d'un faisceau cohérent. Le cas le plus simple est celui où l'on suppose $X$ soit projective, soit de Stein. En effet, si $\fc$ est un faisceau cohérent sur la variété projective ou de Stein $X$, alors il existe une résolution de $\fc$ par des faisceaux localement libres, que l'on peut supposer finie de longueur $n=\dim X$ d'après le théorème des syzygies :
$$0\to E_n\to\cdots\to E_0\to\fc\to 0$$
L'axiome des suites exactes pour les classes de Chern nous amène donc à définir
$$c_{\bullet}(\fc)_{BC,\Z}=c_{\bullet}(E_0)_{BC,\Z}\cdot c_{\bullet}(E_1)_{BC,\Z}^{-1}\cdot \dots \cdot c_{\bullet}(E_n)_{BC,\Z}^{(-1)^n},$$
les classes de Chern étant inversibles dans l'anneau de cohomologie puisque $c_0=1$.

Dans le cas général, on utilise le résultat classique suivant.
\begin{prop} Soit $X$ une variété complexe et $\fc$ un faisceau cohérent sans torsion sur $X$. Alors il existe une modification $\mu:\hat{X}\to X$ telle que $\mu^*\fc$ soit localement libre.\end{prop}
\begin{proof} Au faisceau cohérent $\fc$ on peut associer la variété projective $\Pb(\fc)$ au dessus de $X$, définie comme suit ([Fis76], paragraphe 1.9.) : on a localement une suite exacte
$$\oc_U^{\oplus m}\stackrel{A}{\to}\oc_U^{\oplus n}\to\fc|_U\to 0$$
et la suite duale
$$U\times\C^m\stackrel{^t\!A}{\leftarrow}U\times\C^n\leftarrow V(\fc|_U)\leftarrow 0$$
où les espaces $V(\fc|_U)\subset U\times\C^n$ se recollent en un espace linéaire au dessus de $X$, ce qui permet de définir $\Pb(\fc)\to X$, la variété projective au dessus de $X$ définie par $\fc$.

\end{proof}

On peut de plus supposer $\hat{X}$ lisse à l'aide du théorème d'Hironaka. On peut ainsi définir les classes de Chern d'un faisceau sans torsion, le cas de la torsion étant traité par récurrence.

\pagebreak


\begin{thebibliography}{2}
\bibitem[Aep62]{Aeppli} A. Aeppli, Some exact sequences in cohomology theory for Kähler manifolds, Pacific J. Math., 12, 1962
\bibitem[Dem93]{DemBook} J.-P. Demailly, Analytic Geometry, 1993
\bibitem[EV88]{EV} H. Esnault, E. Viehweg, Deligne-Beilinson cohomology, in : Beilinson conjectures on special values of L-functions, Perspectives in Mathematics 4, Academic press, 1988
\bibitem[Fis76]{Fischer} G. Fischer, Complex analytic geometry, Lecture Notes in Math., 538, Springer, Berlin, 1976
\bibitem[FG86]{Iwasawa} M. Fernandez, A. Gray, The Iwasawa manifold, in : Differential geometry, Peniscola 1985, Lecture Notes in Math., 1909, Springer, Berlin, 1986
\bibitem[Gro58]{Groth} A. Grothendieck, La théorie des classes de Chern, Bull. Soc. math. France, 86, 1958
\bibitem[MK71]{MK} J. Morrow, K. Kodaira, Complex manifolds, Rinehart and Winston, 1971
\end{thebibliography}
\end{document}